\documentclass[12pt]{amsart}
\usepackage{amssymb}
\usepackage{relsize}

\newcommand{\res}{{\mtr{Res}}}

\newcommand{\ovl}{\overline}

\newcommand{\dimccj}{\dim(\cc^j)}

\newcommand{\gr}{\mtr{Gr}}

\newcommand{\grcc}{\gr(\cc)}

\newcommand\numeq[1]%
  {\stackrel{\scriptscriptstyle(\mkern-1.5mu#1\mkern-1.5mu)}{=}}

\newcounter{relctr} 
\everydisplay\expandafter{\the\everydisplay\setcounter{relctr}{0}} 

\AtBeginDocument{} 

\newcommand{\io}{{i_1}}
\newcommand{\itw}{{i_2}}
\newcommand{\jo}{{j_1}}
\newcommand{\jtw}{{j_2}}

\newcommand{\chitw}{\chi_{\itw}}
 
 \newcommand{\chio}{\ch_{\io}}

\usepackage{tikz-cd}
\usepackage{latexsym,amssymb,amsmath}
\usepackage{hyperref}
\usepackage{comment}
\usepackage{mathrsfs}
\usepackage{amsmath,amscd}
\usepackage[enableskew]{youngtab}
\usepackage[all]{xy}
\usepackage{empheq}

\usepackage{amsmath}

\usepackage{empheq}
\usepackage{xcolor}
\definecolor{lightgreen}{HTML}{90EE90}

\newcommand{\cocc}{\co(\cc)}

 \newcommand{\vsk}{\vskip 0.15cm \noindent}

\newcommand{\ccj}{\cc_j}

 \newcommand{\fpcc}{\fp(\cc)}

\newtheorem{theorem}{Theorem}

\newcommand{\sent}{\mapsto}

\newcommand{\mui}{\mu_i}\newcommand{\muj}{\mu_j}

\newcommand\C{\mathcal{C}}
\DeclareMathOperator{\id}{id}

\DeclareMathOperator{\ev}{ev}
\DeclareMathOperator{\coev}{coev}

\excludecomment{verlong}
\includecomment{vershort}
\excludecomment{noncompile}

\usepackage{tikz}
\usetikzlibrary{matrix}
\usepackage{combelow}

\usepackage{amsmath,amsthm,amsfonts,amssymb}
\usepackage{amsmath, amsthm, amssymb, amscd, amsfonts}
\usepackage[all]{xy}
\usepackage{amssymb, amsthm, amsmath}
\usepackage{amsfonts}
\usepackage{color}

\newcommand{\ra}{\rightarrow}

\newcommand{\ot}{\otimes}
\newcommand{\co}{\mathcal O}
\newcommand{\xra}{\xrightarrow}
\newcommand{\mtc}{\mathcal}

\newcommand{\lam}{\lambda}

\newcommand{\Lam}{\Lambda}

\newcommand{\al}{\alpha}
\newcommand{\eps}{\epsilon}

\newcommand{\ul}{\underline}

\newcommand{\lh}{\leftharpoonup}

\newtheorem{lem}[theorem]{Lemma}
\theoremstyle{plain}

\newtheorem{thm}[theorem]{Theorem}
\newtheorem{prop}[theorem]{Proposition}
\newtheorem{defn}[theorem]{Definition}
\newtheorem{cor}[theorem]{Corollary}
\newtheorem{rem}[theorem]{Remark}
\newcommand{\ch}{\chi}
\newcommand{\mtr}{\mathrm}

\newcommand{\ncm}{\newcommand}
\ncm{\np}{\newpage}
\ncm{\ebl}{\end{thebibliography}}
\ncm{\bbl}{\begin{thebibliography}}
\ncm{\chd}{_{ _{\ch}}}
\ncm{\ald}{_{ _{\al}}}
\newcommand{\blam}{\Lam}
\ncm{\cP}{\mathcal{P}}
\ncm{\ei}{e_i}
\ncm{\eij}{e_{i,\;j}}
\ncm{\bt}{\begin{thm}}
\ncm{\bdef}{\begin{defn}}
\ncm{\edf}{\end{defn}}
\ncm{\et}{\end{thm}}
\ncm{\bc}{\begin{cor}}
\ncm{\bl}{\begin{lem}}
\ncm{\el}{\end{lem}}
\ncm{\bpf}{\begin{proof}}
\ncm{\epf}{\end{proof}}
\ncm{\ec}{\end{cor}}
\ncm{\ord}{\mtr{ord}}
\ncm{\er}{\end{rem}}
\ncm{\br}{\begin{rem}}
\ncm{\bn}{\begin}

\ncm{\bp}{\begin{prop}}
\ncm{\ep}{\end{prop}}
\ncm{\bd}{
\begin{document}}
\ncm{\ed}{\end{document}}
\ncm{\beq}{\begin{equation}}
\ncm{\beqn}{\begin{equation*}}
\ncm{\eeq}{\end{equation}}
\ncm{\eeqn}{\end{equation*}}
\ncm{\bea}{\begin{eqnarray}}
\ncm{\eea}{\end{eqnarray}}
\ncm{\beanon}{\begin{eqnarray*}}
\ncm{\eeanon}{\end{eqnarray*}}\ncm{\ek}{\eps|_K}\ncm{\diez}{\#}
\ncm{\bwt}{\bowtie}
\ncm{\cC}{\mtc{C}}\ncm{\cc}{\mtc{C}}
\ncm{\cX}{\mtc{X}}
\ncm{\wt}{\widetilde}
\ncm{\sg}{\sigma}
\ncm{\Rep}{\mathrm{Rep}}
\DeclareMathOperator{\Irr}{Irr}
\ncm{\X}{\mathcal{X}}
\ncm{\cA}{\mathcal{A}}
\ncm{\HKer}{\mtr{HKer}}
\ncm{\LKER}{\mtr{LKer}}
\ncm{\aad}{\mtr{ad}}
\newcommand{\mbf}{\mathbb F}
\ncm{\Dr}{\mtr{D}}
\ncm{\cD}{{\mathcal{D}}}\ncm{\cd}{{\mathcal{D}}}\ncm{\ce}{{\mathcal{E}}}
\ncm{\G}{\mathcal{G}}
\ncm{\Dc}{\mtc{D}}
\ncm{\E}{\mtc{E}}
\ncm{\fp}{\mtr{FPdim}}
\ncm{\Vc}{\mtr{Vec}}
\ncm{\cK}{\mtc{K}}
\ncm{\cM}{\mtc{M}}
\ncm{\cE}{\mtc{E}}
\ncm{\cS}{\mtc{S}}

\newcommand{{\ipr}}{i'}

\DeclareMathOperator{\End}{End}
\ncm{\cop}{\mtr{cop}}
\ncm{\op}{\mtr{op}}
\ncm{\chr}{character }\ncm{\ck}{\mtc{K}}
\ncm{\bw}{\bwt}
\ncm{\hker}{\mtr{HKer}}
\ncm{\bx}{\boxtimes}
\ncm{\blue}{\textcolor[rgb]{.00, .00, 1.00}}
\ncm{\red}{\textcolor[rgb]{1.00, .00, .00}}
\ncm{\green}{\textcolor[rgb]{.50, 0.20, .90}}
\ncm{\bne}{\begin{enumerate}}
\ncm{\ene}{\end{enumerate}}
\ncm{\lker}{\mtr{LKer}}
\ncm{\md}{\medbreak}
\ncm{\rep}{\Rep}\ncm{\ind}{\mtr{ind}}
\ncm{\mdn}{\md\noindent}
\ncm{\dd}{$}
\ncm{\up}{^}
\newcommand{\tcs}{\text}
\newcommand{\mbb}{\mathbb B}
\newcommand{\vs}{\mathbb V}
\newcommand{\sth}{suppose that\;}
\newcommand\rad{\operatorname{rad}}
\newcommand{\itm}{\item}
\newcommand{\dbd}{$$}
\newcommand{\mol}{\mtr{mod}}
 \newcommand{\ro}{\rho}
\newcommand{\irr}{\mathrm{Irr}}
\newcommand{\mbc}{\mathbb C}
\newcommand{\mbs}{\mathbb S}
\newcommand{\mbz}{\mathbb Z}
\newcommand{\ct}{\mtc T}
\newcommand{\sm}{\setminus}
\newcommand{\epl}{^{+}}
\newcommand{\sbsq}{\subseteq}
\newcommand{\sbs}{\subset}
\newcommand{\cco}{\mtr{co}}
\newcommand{\cz}{\mathcal{Z}}
\newcommand{\dual}{^{*}}
\newcommand{\Gm}{\Gamma}
\ncm{\cY}{\mtc{Y}}
\newcommand\ZZ{{\mathbb Z}} 
\newcommand{\bab}{\color{DarkOrchid}{}}
\newcommand{\eab}{\normalcolor{}}
\newcommand{\subs}{\subsection}
\newcommand{\cv}{\mtc{V}}
  \newcommand{\grn}{\green}
\newcommand{\dt}{\delta}

\newcommand{\ccf}{\mathrm{ {CF}(\cc)}}
\newcommand{\cce}{\mathrm{ {CE}(\cc)}}
\newcommand{\cecc}{\mathrm{ {CE}(\cc)}}
\newcommand{\cecd}{\mathrm{ {CE}(\cd)}}
\newcommand{\kk}{\Bbbk}
\newcommand{\otL}{\ot_{L}}
\newcommand{\otl}{\ot_{L}}
\newcommand{\unpsi}{1_{\psi}}
\newcommand{\epsi}{e_{\psi}}
\newcommand{\ephi}{e_{\phi}}
\newcommand{\ech}{e_{\ch}}
\newcommand{\nleftcid}{\text{left normal  coideal subalgebra}}
\newcommand{\dimL}{\dim_{\kk}L}
\newcommand{\cl}{\mtc L}
\newcommand{\mj}{\mtc J}
\newcommand{\tl}{\tilde L}
\newcommand{\tL}{\tilde L}
\newcommand{\tpsi}{\tilde(\psi)}
\newcommand{\tmx}{\tilde{\mtc X}}
\newcommand{\zlh}{\mathrm{ZL}}
\newcommand{\ba}{\mathrm A}
\newcommand{\bv}{\mathrm V}
\newcommand{\zhopf}{\mtc{Z}_{\mtr{Hopf}}}
\newcommand{\lstar}{L^{*}}
\newcommand{\ldstar}{L^{**}}
\newcommand{\mstar}{M^{*}}
\newcommand{\mdstar}{M^{**}}
\newcommand{\lkera}{\lker_{A}}
\newcommand{\mdprime}{M''}
\newcommand{\ldprime}{L''}
\newcommand{\cm}{\mtc M}
\newcommand{\ccm}{\mathcal M}
\newcommand{\cn}{\mathcal N}
\newcommand{\ccn}{\mathcal N}
\newcommand{\rx}{\mtr{Rex}}
\newcommand{\cca}{\ca}
\newcommand{\ih}{\underline{\mtr{Hom}}}
\newcommand{\cih}{\underline{\mtr{coHom}}}
\newcommand{\hm}{\mtr{ {Hom}}}
\newcommand{\cov}{\mtr{coev}}
\newcommand{\rora}{\rho^{\mtr{ra}}}
\newcommand{\rola}{\rho^{\mtr{la}}}
\newcommand{\cx}{\mtc X}
 \newcommand{\cZ}{\cz}
 \newcommand{\ca}{\cA}
 \newcommand{\stat}{\noindent}
 \newcommand{\bfa}{{\bf A}}
 \newcommand{\unu}{\mathbf{1}}
 \newcommand{\barzu}{{\bar {  Z}(\unu)}}
 
\newcommand{\idx}{\id_X}
\newcommand{\lprime}{L'}
\newcommand{\mprime}{M'}
\newcommand{\nat}{ \mtr{{  Nat}}}
\newcommand{\ft}{\mtc F_\lam}
\newcommand{\rhau}{\rightharpoonup}
\newcommand{\lhau}{\leftharpoonup}
\newcommand{\cf}{\mathrm{ {CF}}}

\newcommand{\cfc}{\mathrm{{CF}}(\cc)}
\newcommand{\csu}{\overline{\mathfrak{  C}}}
\newcommand{\cfcc}{\mathrm{ {CF}}(\cc)}
\newcommand{\cfcd}{\mathrm{CF}(\cd)}
\newcommand{\cfd}{\mathrm{CF}(\cd)}
\newcommand{\czcc}{{\cz(\cc)}}
\newcommand{\czcd}{{\cz(\cd)}}
\newcommand{\czt}{{\cz(\cz(\cc))}}
\newcommand{\enx}{\mtr{  End}}
\newcommand{\runu}{R(\unu)}

\newcommand{\bdfn}{\bn{defn}}
\newcommand{\edfn}{\end{defn}}
\newcommand{\deltax}{\delta_X}
\newcommand{\deltav}{\delta_V}
\newcommand{\repcca}{\rep_\cc(A)}
\newcommand{\xotay}{X \ot_A Y}
\newcommand{\xoty}{X \ot Y}
\newcommand{\votw}{V \ot W}
\newcommand{\votaw}{V \ot_A W}
\newcommand{\dimax}{\dim_AX}
\newcommand{\dimccx}{\dim_\cc(X)}
\newcommand{\dimcca}{\dim_\cc(A)}
\newcommand{\dimccv}{\dim_\cc(V)}
\newcommand{\dima}{\dim_A}
\newcommand{\biga}{A}
\newcommand{\comp}{\mathbb C}
\newcommand{\tehtaa}{\theta_A}
\newcommand{\tetaa}{\theta_A}
\newcommand{\ida}{\id_A}
\newcommand{\hma}{\hm_A}
\newcommand{\hmcc}{\hm_\cc}
\newcommand{\fv}{F(V)}
\newcommand{\fw}{F(W)}
\newcommand{\ota}{\ot_A}
\newcommand{\repza}{\rep_\cc^0(A)}
\newcommand{\epsa}{\eps_A}
\newcommand{\bndefn}{\bn{defn}}
\newcommand{\edefn}{\end{defn}}
\newcommand{\bdefn}{\bn{defn}}

\newcommand{\vld}{V^{*}}
\newcommand{\vldd}{V^{**}}
\newcommand{\xld}{X^{*}}
\newcommand{\xldd}{X^{**}}
\newcommand{\yld}{Y^{*}}
\newcommand{\yldd}{Y^{**}}
\newcommand{\aldu}{A^{*}}
\newcommand{\aldd}{A^{**}}

\newcommand{\ia}{\mtr{i}_A}
\newcommand{\aota}{A\ot A}

\newcommand{\idv}{\id_V}

\newcommand{\ld}{^*}
\newcommand{\repg}{\rep(G)}

\newcommand{\thetav}{\theta_V}

\newcommand{\tta}{\theta_A}

\newcommand{\muv}{\mu_V}
\newcommand{\muw}{\mu_W}

\newcommand{\dimcc}{\dim(\cc)}
\newcommand{\chii}{\chi_i}
\newcommand{\chistar}{\ch_{i^*}}
\newcommand{\chj}{\ch_j}
\newcommand{\chm}{\ch_m}
\newcommand{\chn}{\ch_n}
\newcommand{\dimvi}{\dim(V_i)}
\newcommand{\mtcd}{Q}
\newcommand{\mtca}{\mtc A}
\newcommand{\lamcd}{\lam_\cd}
\newcommand{\fpdimcd}{\fp(\cd)}
\newcommand{\laml}{\lam_L}
\newcommand{\apm}{A//M}
\newcommand{\apl}{A//L}
\newcommand{\repapm}{\rep(\apm)}
\newcommand{\repapl}{\rep(\apl)}
\newcommand{\dimvj}{\dim(V_j)}
\newcommand{\dvi}{\dim(V_i)}
\newcommand{\dvj}{\dim(V_j)}
\newcommand{\sumjtom}{\sum_{j=0}^m}
\newcommand{\sumitom}{\sum_{i=0}^m}
\newcommand{\sij}{s_{ij}}
\newcommand{\sji}{s_{ji}}
\newcommand{\dxj}{d_j}
\newcommand{\dxi}{d_i}
\newcommand{\dimka}{\dim_{\kk}(A)}
\newcommand{\dimk}{\dim_{\kk}}
\newcommand{\blaml}{\blam_L}
\newcommand{\sumjtor}{\sum_{j=0}^r}
\newcommand{\dimkl}{\dim_{\kk}(L)}
\newcommand{\mtcjl}{\mtc J_L}
\newcommand{\vota}{ V\ot A}
\newcommand{\vi}{V_i}
\newcommand{\vj}{V_j}
\newcommand{\dimcd}{\dim(\cd)}

\newcommand{\alij}{\al_{ij}}
\newcommand{\alji}{\al_{ji}}
\newcommand{\rcc}{r_\cc}
\newcommand{\rcd}{r_\cd}
\newcommand{\clsx}{[X]}
\newcommand{\clsy}{[Y]}
\newcommand{\clsz}{[Z]}
\newcommand{\rcdp}{r_{\cd'}}
\newcommand{\sumjtorp}{\sum_{j=0}^{r'}}
\newcommand{\aljm}{\al_{jm}}
\newcommand{\aljn}{\al_{jn}}
\newcommand{\sjm}{s_{jm}}
\newcommand{\smj}{s_{mj}}
\newcommand{\snj}{s_{nj}}

\newcommand{\betaij}{\beta_{ij}}
\newcommand{\betaji}{\beta_{ji}}

 \newcommand{\ip}{i'}
\newcommand{\sumjtoprp}{\sum_{j=0}^{r'}}
\newcommand{\sumjtopr}{\sum_{j=0}^{r}}
 \newcommand{\teh}{\tilde{h}}
\newcommand{\cdp}{{\cd'}}\newcommand{\xphii}{X_{\phi(i)}}
\newcommand{\inv}{^{-1}}

\newcommand{\fq}{f_Q}
\newcommand{\tr}{\mtr{tr}}
\newcommand{\rtwone}{R_{21}R}

\newcommand{\ccad}{{\cc_{\mtr{ad}}}}
\newcommand{\ccpt}{{\cc_{\mtr{pt}}}}
\newcommand{\qtr}{quasi-triangular\;}
\newcommand{\trq}{\tr_q}

\newcommand{\repal}{\mtr{Rep}(A//L)}
\newcommand{\lkeravi}{\lker_A(V_i)}
\newcommand{\lkeravj}{\lker_A(V_j)}
\newcommand{\cross}[1][1pt]{\ooalign{%
 \rule[1ex]{1ex}{#1}\cr
 \hss\rule{#1}{.7em}\hss\cr}}
\newcommand{\blml}{\blam_L} 
\newcommand{\phir}{\phi_R}
\newcommand{\kda}{{  \Phi(A)}}

\newcommand{\mtcil}{\mtc{I}_L}

\newcommand{\un}{\unu}
\newcommand{\tfl}{\mtc{T}}
\newcommand{\barzm}{\barz(M)}
\newcommand{\barzn}{\barz(N)}
\newcommand{\ccr}{\mtc R^{\cc}}
\newcommand{\ulc}{\ul{\cc}}

\newcommand{\pimx}{\pi_{M;\;X}}
\newcommand{\pinx}{\pi_{N;\;X}}
\newcommand{\acc}{{\mathrm A_\cc}}
\newcommand{\epsu}{\eps_\unu}

\newcommand{\ob}{\mtr{Obj}}
\newcommand{\obc}{\mtr{Obj(\cc)}}
\newcommand{\ccop}{\cc^{\mtr{op}}}
\newcommand{\mtf}{\mtc F_\lam}
\newcommand{\mtfi}{\mtc F^{-1}_\lam}
\newcommand{\elcd}{\ell_\cd}
\newcommand{\mcid}{\mtc I_\cd}
\newcommand{\mcidp}{\mtc I_{\cd'}}
\newcommand{\wtildelcd}{\widetilde{\elcd}}
\newcommand{\wtildelcdp}{\widetilde{\ell_{\cd'}}}
\newcommand{\cpt}{\cc_{\mtr{pt}}}
\newcommand{\barzr}{\barz_\cd}
\newcommand{\barzv}{\barz(V)}
\newcommand{\acd}{\mathrm A_\cd}
\newcommand{\czrcd}{\cz_\cc(\cd)}
\newcommand{\sml}{\Small}
\newcommand{\bs}{\blue{\Small }}
\newcommand{\yd}{Yetter-Drinfeld}

\newcommand{\sumitor}{\sum_{i=0}^r}
\newcommand{\cdop}{\cd^{\mtr{op}}}
\newcommand{\ccrev}{\cc^{\mtr{rev}}}
\newcommand{\barz}{{\bar{\mathrm Z}}}
\newcommand{\etl}{etale\;}
\newcommand{\czca}{\cz(\ca)}

\newcommand{\tetx}{\text}
\newcommand{\widehta}{\widehat}
\newcommand{\wdhat}{\widehat}
\newcommand{\wht}{\widehat}
\newcommand{\cofa}{{\mathbb C[A]}}
\newcommand{\wdt}{\widehat}
\newcommand{\dl}{{^\#}}
\newcommand{\comx}{\mathbb C}

\newcommand{\sgj}{\sg(j)}

\newcommand{\mujo}{\mu_\jo}
\newcommand{\mujtw}{\mu_\jtw}
\newcommand{\adz}{a^{\#}}
\newcommand{\bdz}{b^{\#}}

\newcommand{\spr}{S^\perp}
\newcommand{\cofs}{\comp [S]}
\newcommand{\spz}{S^{\perp_z}}

\newcommand{\omz}{\omega_z}
\newcommand{\zg}{\mathrm{Z}(S)}
\newcommand{\aling}{{\al \in g}}

\newcommand{\blkg}{\mtr{Bl}(g)}
\newcommand{\clsg}{\mtr{Cl}(g)}
\newcommand{\mtadinv}{\mtc G^{{-1}}}
\newcommand{\muk}{\mu_{k}}
\newcommand{\mta}{\mtc F}
\newcommand{\cofad}{\comp[\wdht A]}
\newcommand{\wtau}{\wdht{\tau}}
\newcommand{\mtainv}{{\mta}^{-1}}
\newcommand{\wdht}{\widehat}
\newcommand{\augm}{\mtr{aug}}
\newcommand{\mua}{\wdht {\wdht a}}
\newcommand{\aps}{A//S}
\newcommand{\cssa}{\cc(S, A)}
\newcommand{\aug}{\mtr{aug}}
\newcommand{\rss}{{\big|_S}}
\newcommand{\gprp}{g^\perp}
\newcommand{\alins}{{s \in S}}

\newcommand{\sz}{s^{D}}
\newcommand{\wmu}{\widehta{\mu}}
\newcommand{\wmui}{\widehta{\mu}_i}
\newcommand{\wmuj}{\widehta{\mu}_j}
\newcommand{\wch}{\widehta{\ch}}

\newcommand{\wzd}{\widehat{d}}
\newcommand{\wpm}{\widehat{P}}
\newcommand{\wps}{\widehat{p}}

\newcommand{\gal}{\mtr{Gal}}
\newcommand{\galkq}{\gal(\mathbb K/\mathbb Q)}
\newcommand{\sgf}{\sg_{_F}}
\newcommand{\sggi}{{\sg(i)}}
\newcommand{\sge}{\sg_{_E}}
\newcommand{\unue}{{\unu_{\cecc}}}

\newcommand{\mtcf}{\mtc {F}}

\newcommand{\wsgf}{\widehat{{\sg}_{F}}}

\newcommand{\we}{\widehta{E}}
\newcommand{\sumktom}{\sum_{k=0}^m}

\newcommand{\wf}{\widehat{F}}

\newcommand{\hsgj}{\widehat{\sg}(j)}
\newcommand{\whsgi}{\widehta{\sg}(i)}

\newcommand{\wpp}{\widehat{p}}
\newcommand{\tauj}{{\tau(j)}}
\newcommand{\dimcctauj}{\dim(\cc^\tauj)}
\newcommand{\etas}{{\eta(s)}}
\newcommand{\mcc}{m_\cc}

\newcommand{\wal}{\widehta{\al}}
\newcommand{\wj}{\widehat{\mtc J}}
\newcommand{\galc}{\mtr{Gal}_{\cc}}
\newcommand{\galz}{\mtr{Gal}_{\czcc}}
\newcommand{\wjr}{\widehat{J}_{R}}

\newcommand{\dimcck}{\dim(\cc^k)}

\newcommand{\wgrcc}{\widehat{\mtr{Gr}(\cc)}}
\newcommand{\nchi}{{\frac{\ch_i}{d_i}}} \newcommand{\nchj}{{\frac{\ch_j}{d_j}}}
\newcommand{\wni}{\widehat{n_i}}

\newcommand{\sgte}{\widetilde{\sg_E}}

\newcommand{\mtad}{\mtc G}
\newcommand{\whj}{\widehta{h}_j}
\newcommand{\jdl}{{j\dl}}
\newcommand{\wcfcc}{\widehat{\cfcc}}
\newcommand{\mutauj}{\mu_{\tau(j)}}
\newcommand{\tauk}{\tau(k)}
\newcommand{\muzm}{{\mu_0^{-}}}
\newcommand{\sqrtog}{\sqrt{|G|}}
\newcommand{\muz}{\mu_0}
\newcommand{\njtw}{n_\jtw}
\newcommand{\njo}{n_\jo}
\newcommand{\fjo}{F_\jo}
\newcommand{\fjtw}{F_\jtw}
\newcommand{\wta}{\widehat{A}}

\newcommand{\dol}{{^{\circ}}}
\newcommand{\bdl}{{b\dl}}
\newcommand{\jdol}{{j\dol}}
\newcommand{\fj}{F_j}

\newcommand{\cwta}{\comp[\wta]}

\newcommand{\hx}{\widehta{x}}
\newcommand{\hy}{\widehta{y}}

\newcommand{\cal}{\mtc A_{\al}}
\newcommand{\xuu}{x_{uu}}
\newcommand{\wxuu}{\widehat{\xuu}}
\newcommand{\xvv}{x_{vv}}
\newcommand{\xuv}{x_{uv}}
\newcommand{\xmn}{x_{m,n}}
\newcommand{\buvmn}{B^{u,v}_{m,n}}
\newcommand{\blm}{\blam}
\newcommand{\dimccr}{\dim(\cc^r)}
\newcommand{\adl}{a\dl}
\newcommand{\sumltom}{\sum_{l=0}^m}

\newcommand{\mbq}{\mathbb Q}
\newcommand{\mbqs}{\mathbb Q(S)}
\newcommand{\mbk}{\mathbb K}
\newcommand{\mz}{\mathbb Z}
\bd
\begin{abstract}
In this paper we show that integral fusion categories with rational structure constants admit a natural group of symmetries given by the Galois group of their character tables. We also generalize a well known result of Burnside from representation theory of finite groups. More precisely, we show that any row corresponding to a non invertible object in the character table of a weakly integral fusion category contains a zero entry.
\end{abstract}
\title[Fusion categories]{On the Galois symmetries for the character table of an integral fusion category}
\author{Sebastian Burciu}
\address{Inst.\ of Math.\ ``Simion Stoilow" of the Romanian Academy P.O. Box 1-764, RO-014700, Bucharest, Romania}
\email{sebastian.burciu@imar.ro}
\date{\today}

\maketitle

\section{Introduction}
Fusion categories can be regarded as a natural generalization of the category of representations  finite groups. From this point of view, many results from group representations have been extended to the settings of fusion categories. For example, in \cite[Theorem 1.6]{eno-weakly} the author showed that any fusion category  of Frobenius-Perron dimension $p^aq^b$ is semi-solvable extending the famous Burnside theorem for finite groups. 

Recently, Shimizu  introduced  in \cite{scalg}  the notion of conjugacy classes for fusion categories, extending  Zhu's work from \cite{zind} and some results for semisimple Hopf algebras obtained previously by  Cohen and Westreich  in  \cite{CW6}.  In the same paper \cite{scalg} the author associated to each conjugacy class a central element called also conjugacy class sum.  These class sums play the role of the sum of group elements of a conjugacy class in group theory and allow one to define character  tables for pivotal fusion categories. In the same paper \cite{scalg}, Shimizu proved the orthogonality relations for the character table, generalizing the famous orthogonality representations for finite groups.

The goal of this paper  is to investigate some other properties of  the character tables of finite groups that can be extended  to the settings of fusion categories.

Let $\cc$ be a pivotal fusion category. Shimizu has constructed in \cite{scalg} a ring of class function $\cfcc$ analogues to the character ring of a semisimple Hopf algebra. He also associated to any object $X$ of $\cc$ a class function $\mtr{ch}(X)\in \cfcc$. More details are given in Section \ref{prelim}. If $X_0, X_1, \dots , X_m$ are a complete set of representatives of the simple objects of $\cc$ we denote by $\ch_i:=\mtr{ch}(X_i)$.

Let $F_0, F_1, \dots ,F_m$ be the central primitive idempotents of $\cfcc$.  Recall that  $\mtr C_j:={\mtf}^{-1}(F_j)\in \cecc$ is called the {\it conjugacy class sums} corresponding to the  conjugacy class $\mtc C_j$ corresponding to $F_j$, see Section \ref{prelim}. Here $\mtf:\cecc\ra \cfcc$ is the usual Fourier transform associated to $\cc$, see  also Section \ref{prelim} for more details.

Since $\{\mtr C_i\}_i$ form a basis for the set of central elements $\cecc$ it follows that the product of two such class sums is a linear combination of the class sums.
Therefore one can write that
$$
\mtr C_i\mtr C_j =\sumktom {c}^{\;k}_{ij} \mtr C_k
$$
where $c^k_{ij}\in \mathbb C$ are complex numbers. It is well-known that for finite groups these structure constants are in fact integer numbers.

Cohen and  Westreich in \cite[Theorem 2.6]{CW6} proved that  for a semisimple quasi-triangular Hopf algebra these structure constants are rational numbers (integers up to a factor of $\dim(H)^{-2}$). Recently in \cite[Theorem 5.6]{zz} the authors improved the above result by showing for a semisimple quasi-triangular Hopf algebra $H$  that this factor can be replaced by $\dim(H)$. Moreover, this result was generalized for premodular fusion categories in \cite{scbf8}. We should mention that at this stage we are not aware of any integral fusion category whose structure constants are not rational numbers.

For any irreducible character $\ch_i$ one  can write $\ch_i=\sumjtom \alij F_j$. Let $\mbk$ the field obtained by adjoining $\alij$ to $\mbq$. By \cite[Corollary  8.53]{eno-annals}  there is a cyclotomic field $ \mathbb Q(\xi)$ such that $\mathbb Q\subseteq \mathbb K\subseteq \mathbb Q(\xi)$.

Our first main result is concerned with the action of the Galois group $\galkq$ on the set of irreducible characters of $\cc$. The notations are explained in details in Section \ref{galois}.

\bt\label{int-main2}
Suppose that $\cc$ is an integral fusion category with rational structure  constants $c^k_{ij}\in \mbq$. Then $\sge$ is an algebra map  for all $\sg \in \galkq$. Moreover, there is a permutation $\eta=\eta_\sg$ such that 
$$
\sgf(\ch_i)=\ch_{\eta(i)},\;\;\sge(E_i)=E_{\eta^{-1}(i)}
$$ 
for all $i \in \mtc I$.
\et
Recall that {\it the character table} of a fusion category $\cc$ with a commutative Grothendieck ring is defined as $\ch_i(C_j)$ where $\ch_i$ are the irreducible characters of $\cc$ and $C_j$ are the class sums of $\cc$ mentioned above.

A classical result of Burnside in character theory of finite groups states that for any  irreducible character $\ch$ of a finite group $G$ with $\ch(1)>1$ there is some
$g \in G$ such that $\ch(g)=0$, see \cite[Chapter 21]{bz}. This result was generalized in \cite[Appendix]{gnn} to weakly integral modular categories.

Our second main result is the following generalization of the above result to weakly-integral fusion categories:
\bt\label{burnside-zero} Let $\cc$ be a weakly-integral fusion category.Then any row corresponding to a non-invertible object in its  character table  contains a zero entry.
\et

Shortly, the organization of this paper is the following. In Section \ref{prelim} we recall the basic on fusion categories and the concept of conjugacy classes and class sums introduced in \cite{scalg}. In  Section \ref{galois} we introduce the Galois group associated to a fusion category and investigate when this group produces symmetries on the irreducible characters of $\cc$. In Section \ref{zeros} we prove theorem \ref{burnside-zero}. We also give a short comparison of our approach with the one from \cite{gnn} for weakly integral modular categories.

All representations and fusion categories in this paper are considered over the ground filed  $\comp$ of complex numbers.

\section{Preliminaries}\label{prelim}
For a finite abelian category $\cc$ we denote by ${\gr}(\cc)$ its Grothendieck group and set ${\gr}_\kk(\cc):={\gr}(\cc)\ot_{\mathbb Z}\kk$. It is well known that for a finite tensor category ${\gr}_\kk(\cc)$ is a $\kk$-algebra with $[V]\cdot [W]=[V\ot W]$ for any two objects $V, W$ of $\cc$.

By a fusion category we mean a semisimple finite tensor  category. We refer to \cite{EGNO15} for the basic theory of tensor categories. For a fusion category $\cc$ we denote by $\irr(\cc)$ the set of isomorphism classes of simple objects of $\cc$. 
It is well known that for a fusion category $\gr(\cc)$ is a based unital ring and therefore one can define Frobenius-Perron dimensions $\fp(X)$ for any object $X$ of $\cc$. The Frobenius-Perron dimension of the category $\cc$ is defined as $$\fp(\cc):=\sum_{X\in \irr(\cc)}\fp(X)^2.$$

A fusion category $\cc$ is called {\it weakly integral} if $\fp(\cc)\in \mathbb Z$. Moreover a fusion category is called {\it integral} if $\fp(X)\in \mathbb Z_{>0}$ for any simple object $X$ of $\cc$.

Throughout this paper $\cc$ denotes a fusion category and $\unu$ the unit object of  a $\cc$. The {\it monoidal centre} (or Drinfeld centre) of $\cc$ is a braided fusion category $\czcc$ constructed as follows, see e.g. \cite[XIII.3]{Kas} for details. The objects of $\czcc$ are pairs $(V, \sg_V)$ of an object $V\in \cc$ and a natural isomorphism
$\sg_{V, X}: V \ot  X \ra X \ot V $ for all $X \in \cocc$, satisfying a part of the hexagon axiom. A morphism  $f:(V, \sg_V)\ra (W, \sg_W)$ in $\czcc$
 is a morphism in $\cc$  such that $(id_X\ot f) \circ\sg_{V, X}=\sg_{W, X}\circ(f\ot id_X)$ for all objects $X$ of $ \cc$. The composition of morphisms in $\czcc$ is defined via  the usual composition of morphisms in $\cc$.

Let $\cc$ be a finite tensor category and $F:\czcc\ra \cc$ be the forgetful functor, $F(V, \sg_V)=V$. Then $F$ admits a right adjoint functor $R:\cc \ra \czcc$  such that  $Z :=FR:\cc \ra \cc$ is a Hopf comonad.  Moreover, as in \cite[Subsection 2.6]{scalg} one  has that 
\beq
Z(V)\simeq \int_{X\in \cc}X\ot V\ot X^*.
\eeq
 The counit $\eps:Z\ra \id_\cc $  is given by $\eps_V:= \pi_{V ;1} $ where $\pi_{V;X}:Z (V)\ra X\ot V\ot X^*$ are the universal dinatural transformation associated to the above end $Z (V)$. Moreover, using Fubini's theorem for ends, the comultiplication $\delta : Z  \ra {Z }^2$ of $Z $ is also  described in terms of the dinatural transformation $\pi$ associated to $Z$, see \cite[Subsection 3.2]{scalg}.

The object $A:=Z(\unu)\in \czcc$ has the structure of  a central commutative algebra in $\cz(\cc)$. It is called the {\it adjoint algebra} of $\cc$.

The multiplication and the unit of $A$ are uniquely determined by  by the universal property of the end $Z$ as:
\beqn 
\pi_{\unu;X} \circ m_A = (\id_X \otimes \ev_X \otimes \id_{X^*} ) \circ (\pi_{\unu;X} \otimes \pi_{\unu;X}), \;\;\pi_{\unu;X} \circ u_A = \cov_X.
\eeqn

Recall that a {\it pivotal structure} $j$ on $\cc$ is a tensor isomorphism $j:\id_\cc\ra ()^{**}$.  Given such a pivotal structure one can construct for any object $X$ of $\cc$ a {\it right evaluation}
$\widetilde{ev}_X:X\ot X^*\xra{j\ot id}X^{**}\ot X^*\xra{ev_{X^*}} \unu$. 
For any morphism $f:A\ot X\ra B\ot X$ in $\cc$ one can define  the right partial pivotal trace $ \tr_{A, B}^X\in \hm_\cc(A, B)$ of $f$:
\beqn
\tr_{A, B}^X:A=A\ot \unu\xra{\id \ot coev_X}A\ot X \ot X^*\xra{f \ot id} B\ot X\ot X^*\xra{\\id \ot \widetilde{ev}_X} B.
\eeqn

Then the usual {\it right pivotal trace} of an endomorphism $f :X\ra X$ is obtained as a particular case for $A=B=\unu$. In particular,  the {\it right  (quantum) dimension of $X$} with respect to $j$ is defined as the right trace of the identity of $X$. 

A  pivotal structure (or the underlying fusion category) is called {\it spherical} if
$\dim(X) = \dim(X^*)$ for all objects $X$ of $\cc$, see \cite[Definition 4.7.14]{eno-annals}.

\br\label{dimst}
By \cite[Proposition 29]{eno-annals} one has that $\dim(V^*)=\ovl{\dim(V)}$ for any object of a pivotal fusion category. In particular, in a spherical category $\dim(V )$ is (totally) real number.
\er

Recall also that the global dimension of the pivotal fusion category $\cc$ is defined as $\dimcc:=\sumitom d_id_i^*$. It is well-known that over the complex numbers one has  $\dimcc\neq 0$ for any pivotal fusion category.

A fusion category $\cc$ is called {\it pseudo-unitary} if $\fpcc =\dimcc$. If this is the case,
then by \cite[Proposition 8.23]{eno-annals}, $\cc$ admits a unique spherical structure with respect to which the categorical dimensions of simple objects are all positive. It is called the {\it canonical spherical structure}. For this spherical structure, the categorical dimension of an object coincides with its Frobenius-Perron dimension, i.e. $\fp(X)=\dim(X)$ for any object $X$ of $\cc$.

If $\cc$ is a weakly-integral fusion category then $\cc$ is pseudo-unitary by \cite[Proposition 8.24]{eno-annals}. 


\br\label{wintccjint}
Let $\cc$ be a weakly integral fusion category.  By \cite[Proposition 8.27]{eno-annals} the dimensions of  simple objects in $\cc_{ad}$ are integers. Since $\cc^j$ are sum of simple object of the adjoint subcategory it follows that $\dimccj$ are integers.
\er

The internal character $\mtr{ch}(X)$ of an object $X$ of $\cc$ is defined as the  partial pivotal trace 
\beqn
\mtr{ch}(X):=\tr^{X}_{A, \unu}(\rho_{X}):A\ra \unu.
\eeqn
where $\rho_X:A\ot X\ra X$ is the canonical action of $A$ on $X$, see \cite[Definition 3.3]{scalg} for details.

The {\it space of class functions} of $\cc$ is defined as   $\cfcc:=\hm_\cc(A, \unu)$ and it is a $\comp$-algebra. The multiplication of two class functions $f, g\in \cfcc$ is defined via $f\star g:=f \circ Z(g) \circ \delta_{\unu}.$ Here the map $\delta: Z \ra  Z^2$ is the comultiplication structure of the Hopf comonad $Z$ recalled above.

By \cite[Theorem 3.10]{scalg} for a pivotal fusion category $\cc$ one has that $\mtr{ch}(X\ot Y)=\mtr{ch}(X)\mtr{ch}(Y)$ for any two objects $X$ and $Y$ of $\cc$ and  $\mtr{Gr}_{\comp}(\cc)\ra \cfcc, [X]\sent \mtr{ch}(X)$ is an isomorphism of algebras. Moreover, the character of the unit object $\mtr{ch}(\unu)=\eps_\unu$ is the unit of $\cfcc$. Recall from above that  $\eps:Z\ra \id_\cc $ is the counit of the central Hopf comonad.

By \cite[Theorem 5.9]{scalg} for a pivotal fusion category it follows that the map $\mtr{ch}:\mtr{Gr}_{\kk}(\cc)\ra \cfcc, [X]\sent \ch(X)$ is an isomorphism of $\comp$-algebras.

The  space $\cecc:= \hm_{\C}(\unu, A) $ is called {\it the space of central elements.} It is also a $\comp$-algebra where the multiplication on $\cecc$ is given by
 $a.b :=m \circ (a \ot b)$ for any $a, b \in \cecc$. There is a non-degenerate pairing  $\langle\;,\; \rangle : \cfcc \times \cecc\ra \comp$, given by $ \langle f, a\rangle  \id_{\unu}= f \circ a,$ for all $f \in \cfcc$ and $a\in \cecc$. We also denote $f(a):=\langle f,\;a\rangle$. There is a right action of $\cecc$ on $\cfcc$ denoted by $\lh$  given by $f \lh b=f \circ m \circ (b\ot \id_{A})$ for all $f \in \cfcc$ and $b \in \cecc$.
\subsection{Orthogonality relations for pivotal fusion categories}

For the rest of this paper we fix a pivotal fusion category $\cc$ with commutative Grothendieck ring and we use  following notations. We denote by  $M_0, M_1, \dots ,M_m$   a complete set of representatives for the isomorphism classes of simple objects of $\cc$. For brevity,  we also denote $\ch_i:=\mtr{ch}(M_i)\in \cfcc$ and $d_i:=\ch_i(1_\cecc)$ the categorical dimensions of the simple objects. We call $\ch_i$ the {\it irreducible character} of the simple object $M_i$. 
Without loss of generality we may suppose that $M_0=\unu$ and therefore $d_0=1$. {In this case by \cite[Lemma 6.10]{scalg} one has $\ch_0=\eps_\unu$.} Moreover we denote by $i^*$ the unique index for which $M_i^*\simeq M_{i^*}$.

\br\label{dimnz}
By \cite[Proposition 4.8.4]{EGNO15} it also follows that $d_i\neq 0$ for all $i$, since the dimensions of simple objects are not zero.
\er
Let $\cc$ be a pivotal fusion category with commutative Grothendieck ring. Recall that $R:\cc \ra \czcc$ is the right adjoint of the  forgetful functor $F:\czcc \ra \cc$. Note that in this case by \cite[Theorem 6.6]{scalg} the object $R(\unu) \in \co(\czcc)$ is multiplicity-free.

A {\it conjugacy class} of $\cc$ is defined as a simple subobject of $R(\unu)$ in $\czcc$. Since the monoidal center $\czcc$ is also a fusion category we can write $R(\unu)=\bigoplus_{j=0}^m\mtc C^j$ as a direct sum of simple objects in $\czcc$. Thus  $\mtc C_{0},\dots, \mtc C_{m}$ are the conjugacy classes of $\C$.  Since the unit object $\unu_{\czcc }$ is always a subobject of $R(\unu)$, we may assume  for the rest of this paper that $\mtc C_{0} = \unu_{\czcc }$.


By \cite[Theorem 3.8]{scalg} there is a canonical isomorphism
\beq\label{enx}
\enx_{\czcc}(R(\unu))\simeq \cfcc,\;\; f \sent Z(f)\circ \delta_{\unu}.
\eeq
This isomorphism gives a natural bijection between $F_j$ and $\cc^j$. Then the $\cc^0$ corresponds $F_0=\lam$.

A cointegral $\lam$ of $\cc$ is defined as the unique element (up to a scalar) of $\cfcc$ with the property that $\ch\lam=\ch(\unu)\lam$ for all $\ch \in \cfcc$.

Following \cite{scalg} with the above notations one has
$$\lam=\frac{1}{\dimcc}\big(\sumitom d_i\ch_i\big)\in \cfcc.$$  Moreover, one has $<\lam,u_A>=1$ and $\lam^2=\lam$ where $u_A:\unu \ra A$ is the unit of $A$.

The {\it Fourier transform} of $\cc$  associated to $\lambda$ is the linear map
$
\mtc F_{\lambda}:\cecc\ra \cfcc\;\;\text{given by}\;\;a \mapsto \lambda \lh \mtc S(a)
$
where $\mtc S:\cecc\ra \cecc$ on $\cecc$ is the antipodal map of $\cecc$, see \cite[Definition 3.6]{scalg}.

Let also $F_0, F_1, \dots ,F_m$ be the central primitive idempotents of $\cfcc$. Without loss of generality we may suppose that $F_0=\lam$. We define $\mtr C_j:={\mtf}^{-1}(F_j)\in \cecc$ to be the {\it conjugacy class sums} corresponding to the  conjugacy class $\mtc C_j$. It follows as in \cite[Section 6]{scalg} that $C_0=\mtf^{-1}(\lam)=1_\cecc=u_A$.


Since $\cfcc$ is commutative it follows that the set $\{F_j\}$ of central primitive idempotents forms also a bases for $\cfcc$. Therefore for any irreducible character $\ch_i$ one can write
\beq\label{chi}
\ch_i=\sumjtom \alij F_j
\eeq
for some scalars $\alij \in \comp$.

For $0\leq j\leq m$ let $\muj:\cfc\ra \comp$  be the characters of $\cfcc$ corresponding to the primitive central idempotents $F_j$. Since we assume $\cfcc$ is commutative the characters $\muj$ are  algebra morphisms.

By \cite[Equation (6.13)]{scalg} one has that 
\beq\label{normcj}
\alij=\muj(\ch_i)=\frac{\ch_i(C_j)}{\dim(\cc^j)}.
\eeq

Since in this case $\cc^0=\unu_\czcc$ it follows from the above equation that
\beq\label{mu0}
\mu_0(\ch_i)=\ch_i(1_\cecc)=d_i
\eeq
for all $i$.
Note also that in this case $\mu_0=d$, the quantum dimension morphism.

Following \cite[Example 4.4]{scalg} one has that
$$A=\oplus_{i=0}^m M_i\ot M_i^*$$

In this case the irreducible characters $\ch_i$ are given by
$$
\ch_i:=A\xra{\text{projection}} M_i\ot M_i^*\xra{\tilde{ev}_{M_i}} \unu.
$$

As in \cite[Subsection 6.1]{scalg} one has that the primitive central idempotents $E_i\in \cecc$  of $\cecc$ are given by 
$$
\unu \xra{\coev_{M_i}} M_i\ot M_i^*\hookrightarrow A.
$$ 

Moreover, in this case
$$
<\ch_i, E_j>=\delta_{i,j}
$$
for all $i,j$.

By \cite[Corollary 6.11]{scalg} the first orthogonality relation for $\cc$ can be written as
\beq\label{611}
\sum_{ k=0}^m \dim(\cc^k) \al_{i k} \al_{m^* k}=\delta_{i,m}\dimcc.
\eeq
and the second orthogonality relation as:
\beq\label{612}
\sum_{ i=0}^m\al_{il}\al_{i^*k}=\delta_{l,k}\frac{\dimcc}{\dim(\cc^k)}.
\eeq

Recall from \cite[Equation (3.4)]{scbf8} that
\beq\label{clseqcc}
\sumjtom \frac{\dim(\cc)}{\dim(\cc^j)} F_j\ot F_j=\sumitom \ch_i\ot\ch_{i^*}.
\eeq
From here it follows that 
\beq\label{fj}
F_j:=\frac{1}{n_{j}}\big(\sumitom \ovl{\muj(\ch_{i})}\ch_{i}\big)\eeq 
where
\beq\label{njfus}
n_j=\frac{\dim(\cc)}{\dim(\cc^j)}.
\eeq
For any two class functions $\ch=\sumitom \al_i\ch_i$ and $ \mu=\sumitom \beta_i\ch_i$ of $\cfcc$ we define
$$
m_\cc(\ch, \mu):=\sumitom\al_i\beta_i.
$$
\bl 
Let $\cc$ be a pivotal fusion category. With the above notations it follows that
\beq\label{mtau}
m_\cc(\ch, \mu)=\sumjtom\frac{\dimccj}{\dimcc}\muj(\ch)\ovl{\muj(\mu)}
\eeq
for any two class functions $\ch, \mu\in \cfcc$.
\el
\bpf
Indeed, for any two irreducible characters $\ch=\ch_s$ and $\mu=\ch_t$ the right hand side of the above equation can be written as
$$
\sumjtom\frac{\dimccj}{\dimcc}\muj(\ch_s)\ovl{\muj(\ch_t)}\numeq{\ref{normcj}}\sumjtom\frac{\dimccj}{\dimcc}\al_{sj}\al_{t^*j}\numeq{\ref{611}}\delta_{s,t}.
$$
then the Lemma follows by linearity.
\epf
Note also that $\muj(\ch_{i^*})=\ovl{\muj(\ch_i)}$ for all $i,j$.

In particular, since $\cc^0=\unu_\czcc$ it follows that $n_0=\dimcc$.

Since $\{C_j\}_j$ form a $\comp$-linear basis for $\cecc$ one has that 
\beq\label{defstrc}
\mtr C_\jo \mtr C_\jtw=\sum_{l=0}^m c^l_{\jo,\jtw}\mtr C_l
\eeq
form some scalars $c^l_{\jo,\jtw}\in \comp$. These scalars are called {\it the structure constants of $\cc$.}

By \cite[Theorem 1.1]{scbf8} if $\cc$ is a pivotal fusion category with a commutative Grothendieck ring then
\beq\label{bnf}
c^k_{ij}=\sum_{s=0}^m \frac{\ch_s( \mtr C_i)\ch_s(\mtr C_j)\ch_{s^*}(\mtr C_k)}{\dimcc \dim(\cc^k)d_s}.
\eeq
Using Equation \eqref{normcj} this can be written as
\beq\label{bnf2}
c^k_{ij}=\big(\frac{\dim(\cc_\jo)\dim(\cc_\jtw)}{{\dimcc}}\big)\big(
\sum_{s=0}^m \frac{\al_{si}\al_{sj}\ovl{\al_{sk}}}{d_s}\big).\eeq

\section{Galois action for integral fusion categories}\label{galois}
Let $\cc$ be a pivotal fusion category with a commutative Grothendieck ring $\mtr{Gr}_{\kk}(\cc)\simeq \cfcc$. 

Let also $M_0, M_1, \dots ,M_m$ be a complete set of representatives for the isomorphism classes of simple objects of $\cc$. As above, without loss of generality we may assume that $M_0=\unu$ is the unit object of $\cc$. Moreover, we denote by $\ch_i:=\mtr{ch}(M_i)$ the characters of the simple objects $M_i$. Recall that in this case $\ch_0=\eps_\unu$ is the unit of the algebra $\cfcc$.
 
As before, by $F_0, F_1, \dots ,F_m$ are denoted  the central primitive  idempotents of $\cfcc$ and by $\mu_0, \mu_1, \dots ,\mu_m:\cfcc \ra \comp$  their corresponding characters on $\cfcc$. Therefore $\mui:\cfcc\ra \comp$ are algebras maps and $\mui(F_j)=\delta_{i,j}.$ We also denote by $\cc^0, \cc^1,\dots ,\cc^m$ the conjugacy classes of $\cc$ corresponding in this order to the primitive idempotents $F_0, F_1, \dots ,F_m$. Moreover, we let $\mtr C_0, \mtr C_1,\dots ,\mtr C_m$ denote their corresponding class sums.

Without loss of generality we may suppose that $\mu_0=d$ is the quantum dimension homomorphism. Therefore $\mu_0(\ch_i)=d_i$ for all $i$. It follows as in the previous section that $\mtr C_0=\mtf^{-1}(\lam)=1_\cecc=u_A$, the unit $u_A:\unu\ra A$ of the adjoint algebra $A$.

Writing $\ch_i=\sumjtom\alij F_j$ we let $\mathbb K$ be the field extension of $\mathbb Q$ by all the scalars $\alij$. By \cite[Corollary  8.53]{eno-annals}  there is a cyclotomic field $ \mathbb Q(\xi)$ such that $\mathbb Q\subseteq \mathbb K\subseteq \mathbb Q(\xi)$ fro some root of unity $\xi$.

Define $G:=\galkq$ and note that is an abelian group as a quotient of the abelian group $\mtr{Gal}(\mathbb Q(\xi)/\mbq)$.

\br\label{alijjp}
Note that if $\alij=\al_{ij'}$ for all $i$ then $\mu_j(\ch_i)=\mu_{j'}(\ch_i)$, thus $\mu_j=\mu_{j'}$ and therefore $j=j'$.
\er

\subsection{Action $\wsgf$ on $\wcfcc$}
Recall that $\wcfcc$ is defined as the linear dual space of $\cfcc$. 
For $\muj:\cfcc\ra \comp$ and $\sg \in \galkq$ define  $\sg.\muj:\in \widehat{\cfcc}$ as the linear function on $\cfcc$ which on the basis given by the irreducible characters $\{\ch_i\}_i$ is given by 
\beq\label{def1}
[\sg.\muj](\ch_i)=\sg(\muj(\ch_i))=\sg(\alij).
\eeq
\bl 
For any pivotal fusion category $\cc$ with the above notations it follows that $\sg.\muj:\cfcc\ra \comp$ is an algebra map. 
\el
\bpf
Indeed, suppose that $\chio\chitw=\sumktom N^k_{\io,\itw}\ch_k$. Then one has
$$
[\sg.\muj](\chio\chitw)=[\sg.\muj](\sumktom N^k_{\io,\itw}\ch_k)=\sumktom N^k_{\io,\itw}\sg.\muj(\ch_k)=\sumktom N^k_{\io,\itw}\sg(\muj(\ch_k)).
$$
On the other hand,
\begin{eqnarray*}
[\sg.\muj](\chio)[\sg.\muj](\chitw) &=& \sg(\muj(\chio))\sg(\muj(\chitw))=\sg(\mu_j(\chio)\muj(\chitw)\\ &= & \sg(\muj(\chio\chitw))=\sg(\muj(\sumktom N^k_{\io,\itw}\ch_k))\\ &= &\sumktom N^k_{\io,\itw}\sg(\muj(\ch_k)).
\end{eqnarray*}
\epf
It follows that the group $\galkq$ acts on the set of all algebra homomorphisms $\muj:\cfcc\ra \comp$. It is easy to see  that
$\sg.(\sg'.\muj)=(\sg\sg').\muj.$
Indeed one has 
\begin{eqnarray*}
[\sg.(\sg'.\muj)](\ch_i)& = & \sg\big([\sg'.\muj](\ch_i)\big)=\sg\bigg(\sg'\big(\muj(\ch_i)\big)\bigg)\\ &= &\sg\sg'(\muj(\ch_i))=[(\sg\sg').\muj](\ch_i).
\end{eqnarray*}

We denote by $\mtc J$ the set of all indices $j\in \{0,\dots , m\}$ such that $\muj$ is an algebra map.

Since $\sg.\muj:\cfcc \ra \comp$ is an algebra homomorphism it follows that there is an index $\tau_\sg(j)\in \mtc J$ such that $\sg.\muj=\mu_{\tau_\sg(j)}$. It follows that 
\beq\label{alsg}
\sg(\al_{ij})=\al_{i\tau_\sg(j)}
\eeq
for any $i,j$. Indeed, $\sg(\alij)=\sg.\muj(\ch_i)=\mu_{\tau_\sg(j)}(\ch_i)=\al_{i\tau_\sg(j)}$.

\bp\label{tauperm}
For any $\sg \in \galkq$ one has that 
$\tau_\sg$ is a permutation of $\mtc J$.
\ep 
\bpf
Indeed, if $\tau_\sg(j)=\tau_\sg(j')$ then since $\al_{i\tau_\sg(j)}=\al_{i\tau_\sg(j')}$ it follows that $\sg(\alij)=\sg(\al_{ij'})$ which shows that $\alij =\al_{ij'}$ and therefore  $j=j'$ by Remark \ref{alijjp}.
\epf

Define also a $\comp$-linear map $\wsgf:\wcfcc \ra \wcfcc$ written on the linear $\comp$-basis $\{\mu_j\}_{j\in \mtc J}$ by 
$$
\wsgf(\mu_j)=\sg.\mu_j=\mu_{\tau_{\sg}(j)}.
$$

\br\label{unital}
Note that the map $\wsgf$ preserves the unit $\mu_0$ of $\wcfcc$ if and only if $\sg.\mu_0=\mu_0\iff d_i=\sg(d_i)$ for all $i$.

Then  it follows that all the maps $\wsgf$ are unital (for any $\sg \in \galkq$) if and only if $\cc$ is an integral category.
\er
\bp \label{galact}
With the above notations one has that $\tau_{\sg}\tau_{\sg'}=\tau_{\sg\sg'}$ for all $\sg, \sg'\in \galkq$. This shows that we have a group embedding
$$
\galkq \hookrightarrow \mathbf S_m,\;\; \sg \sent \tau_\sg.
$$
\ep

\bpf
Indeed one has 
$$
\mu_{[\tau_{\sg}\tau_{\sg'}](j)}=\sg(\mu_{\tau_{\sg'}(j)})=\sg.(\sg'.(\muj))=(\sg\sg').(\muj)=\mu_{\tau_{\sg\sg'}(j)}$$
which implies $(\tau_{\sg}\tau_{\sg'})(j)=\tau_{\sg\sg'}(j)$, i.e $\tau_{\sg}\tau_{\sg'}=\tau_{\sg\sg'}$.
\epf
\br \label{tauinv}
The above proposition implies that $\tau_{\sg^{-1}}=\tau^{-1}_{\sg}$, i.e.
\beq\label{alsg2}
\sg^{-1}(\al_{ij})=\al_{i\tau_{\sg}^{-1}(j)}
\eeq
\er
For the rest of the paper we shortly write $\tau:=\tau_\sg$ if $\sg$ is implicitly understood.
\br\label{ncommcfcc}
Note that the permutation action $\tau_{\sg}$ holds also in the case of a non-commutative $\cfcc$ since $\sg.\muj$ is a trace.
\er

\bp \label{sgccj}
Let $\cc$ be a pivotal fusion category with a commutative $\cfcc$.
For any $\sg \in \galkq$ with the above notations one has: 
\beq\label{sgccj}
\sg \bigg( \frac{\dimcc}{\dimcck} \bigg)=\frac{\dimcc}{\dim(\cc^{\tau(k)})}
\eeq
In particular if $\cc$ is  weakly integral one has 
\beq\label{sgccj2}
\dimcck=\dim(\cc^{\tau(k)}).
\eeq
\ep
\bpf 
Since $\cc$ is pivotal, note that the size $\dimccj$ of $\ccj$ is non-zero  since it is the quantum dimension of a simple object in a pivotal fusion category $\czcc$.
By the second orthogonality  relation \cite[Cor. 6.11]{scalg} one has:
\beq\label{orthos2t}
\sum_{ i=0}^m\al_{il}\al_{i^*k}=\delta_{l,k}\frac{\dimcc}{\dimcck}
\eeq
Applying $\sg\in \galkq$ to this Equation one obtains:
\beq\label{orthosg}
\sum_{ i=0}^m\sg(\al_{il})\sg(\al_{i^*k})=\delta_{l,k}\sg(\frac{\dimcc}{\dimcck})
\eeq
This implies that
$$
\sum_{ i=0}^m\al_{i\tau(l)}\al_{i^*\tau(k)}=\delta_{l,k}\sg(\frac{\dimcc}{\dimcck})
$$
On the other hand, by the same orthogonality relation we have:
$$
\sum_{ i=0}^m\al_{i\tau(l)}\al_{i^*\tau(k)}=\delta_{\tau(l), \tau(k)}\frac{\dimcc}{\dim(\cc^{\tau(k)})}
$$
Therefore for $l=k$ it follows that
$$
\sg(\frac{\dimcc}{\dimcck})=\frac{\dimcc}{\dim(\cc^{\tau(k)})}
$$
Moreover, if $\cc$ is weakly integral it follows that 
$$\sg(\dimcck)=\dim(\cc^{\tau(k)}).$$
By Remark \ref{wintccjint}  both dimensions above are integers and therefore $\dimcck=\dim(\cc^{\tau(k)})$.
\epf

\br \label{sgnj}
From this it follows by Equation \eqref{njfus} that
\beq\label{sgnj}
\sg(n_j)=n_{\tau(j)}.
\eeq

In particular, if $\cc$ is weakly integral then $\dim(\cc^j)$ is an integer and it follows that $n_j=n_{\tau(j)}$.
\er
\bp\label{strcint}
If $\cc$ is an integral fusion category then 
\beq\label{bnf3}
\sg(c^k_{ij})=c^{\tau(k)}_{\tau(i)\tau(j)}
\eeq
for all $i,j,k$.
\ep
\bpf
If $\cc$ is integral then $\dim(\cc^j)$ are integers by Remark \ref{wintccjint}. On the other hand $\sg(d_i)=d_i$ since $d_i=\fp(X_i)\in \mathbb Z_{>0}$. Then applying $\sg$ to the above formula it follows that
\begin{eqnarray*}
\sg(c^k_{ij})&=&\big(\frac{\dim(\cc_\jo)\dim(\cc_\jtw)}{{\dimcc}}\big)\big(
\sum_{s=0}^m \frac{\sg(\al_{si})\sg(\al_{sj})\ovl{\sg(\al_{sk})}}{\sg(d_s)}\big)=
\\&\numeq{\ref{alsg}, \ref{sgccj2}}& 
\big(\frac{\dim(\cc_\jo)\dim(\cc_\jtw)}{{\dimcc}}\big)\big(
\sum_{s=0}^m \frac{\al_{s\tau(i)}\al_{s\tau(j)}\ovl{\al_{s\tau(k)}}}{d_s}\big)
\\&=&
c^{\tau(k)}_{\tau(i)\tau(j)}.
\end{eqnarray*}
\epf

\subsection{Action $\sgf$ on $\cfcc$ }
For any $\sg \in \galkq$  we define $\sg(\ch_i)\in \cfcc$ with $\sg.\ch_i:=\sumjtom \sg(\alij)F_j\in \cfcc$ for all irreducible characters $\ch_i$.

Define also a $\comp$-linear map $\sgf:\cfcc \ra \cfcc$ written on the linear $\comp$-basis $\{\ch_i\}_{i}$ by 
$$
\sgf(\ch_i)=\sg.\ch_i=\sumjtom \sg(\al_{ij})F_j,\;\text{if}\;\ch_i=\sumjtom \alij F_j.
$$
\bp \label{sgfonfj}
Let $\cc$ be a pivotal fusion category. With the above notations, for any $\sg \in \galkq$ one has that 
\beq\label{sgfj}
\sgf(F_{j})=F_{\tau^{-1}(j)},\; \sgf^{-1}(F_{j})=F_{\tau(j)},\;\sgf^{-1}=(\sg^{-1})_F.
\eeq

\ep
\bpf 
Recall that $\al_{i\tau(j)}=\mu_{\tau(j)}(\ch_{i})=(\sg.\muj)(\ch_{i})=\sg(\al_{ij})$
and by Equation \eqref{orthos2t}
\beqn
\sum_{ i=0}^m\al_{il}\al_{i^*k}=\delta_{l,k}\frac{\dimcc}{\dimcck}=\delta_{l,k}n_{l}.
\eeqn

On the other hand, since by Equation \eqref{fj} we have the formula $F_j:=\frac{1}{n_{j}}\bigg(\sumitom \ovl{\muj(\ch_{i})}\ch_{i}\bigg)$ it follows that,
\begin{eqnarray*}
\sgf(F_{j}) &=& \frac{1}{n_{j}}\bigg(\sumitom \ovl{\muj(\ch_{i})}\sgf(\ch_{i})\bigg)
\\ &= &\frac{1}{n_{j}}\bigg(\sumitom \ovl{\muj(\ch_{i})}\big(\sumltom\sg(\al_{il})F_{l}\big)\bigg)
\\ &= & 
\frac{1}{n_{j}}\bigg(\sumltom \big(\sumitom\ovl{\muj(\ch_{i})}\sg(\al_{il})\big)F_{l}\bigg)
\\ &= & 
\frac{1}{n_{j}}\bigg(\sumltom \big(\sumitom \ovl{\al_{ij}}\al_{i\tau(l)}\big)F_{l}\bigg)\numeq{\ref{orthos2t}}F_{\tau^{-1}(j)}
\end{eqnarray*}
Since $\tau$ is a permutation of $\mtc J$ it follows that $\sgf$ is bijective. Moreover since $(\tau_{\sg})^{-1}=\tau_{{\sg}^{-1}}$ the last two equality also follow.
\epf

\bp\label{sgfalgm}
Let $\cc$ be a pivotal fusion category. One has that $\sgf:\cfcc\ra \cfcc$ is an algebra map. 
\ep
\bpf
Indeed, suppose that $\chio\chitw=\sumktom N^{k}_{\io \itw}\ch_{k}$. Expanding on the formula $\ch_{i}=\sumjtom\al_{ij}F_{j}$ it follows that
\beq\label{atfj}
\al_{\io j}\al_{\itw j}=\sumktom N^{k}_{\io\itw}\al_{k j}
\eeq
Applying $\sg$ to this Equation one has
\begin{eqnarray*}
\sgf(\chio\chitw) &=& \sumktom N^{k}_{{\io \itw}}\sgf(\ch_{k})
=\sumktom N^{k}_{{\io \itw}}\big(\sumjtom \sg(\al_{kj})F_{j}\big)
\\ &=& \sumjtom\big(\sumktom N^{k}_{{\io \itw}}\sg(\al_{kj})\big)F_{j}
\\ &\numeq{\ref{atfj}}& \sumjtom \sg(\al_{\io j})\sg(\al_{\itw j})F_{j}
\\ &=& \sgf(\ch_{\io})\sgf(\ch_{\itw}).
\end{eqnarray*}
Note that $\ch_0$, the character of the unti object is the identity of the $\mathbb C$-algebra $\cfcc$.

It follows that $\sgf(\ch_{0})=\sgf(\sumjtom F_{j})=\sumjtom F_{\tau^{-1}(j)}=\ch_{0}$ which shows that $\sgf$ is unitary.
\epf
\br
Since $\sgf$ is linear it follows from above  that 
\beq\label{sgfidfj}
\sgf(\sumjtom\al_{j}F_{j})=\sumjtom\al_{j}F_{\tau^{-1}(j)}.
\eeq
for all scalars $\al_{j}\in \comp$.
\er 
\bc 
Let $\cc$ be a pivotal fusion category. For any $\sg, \sg'\in \galkq$ one has 
$$(\sg'\sg)_F=\sg_F\circ \sg'_F.$$
\ec
\bl
If $\cc$ is a weakly integral fusion category then with the above notations one has that:
\beq\label{selfm}
m_\cc(\sg(\ch_\io), \sg(\ch_\itw))=\delta_{\io,\itw}
\eeq
for all $i$.
\el
\bpf 
Since $\cc$ is weakly integral one has that $\dim(\cc^j)=\dim(\cc^{\tau(j)})$ for all $j$.
On the other hand  by Equation \eqref{mtau} one has that
\begin{eqnarray*}
\mcc(\sg(\ch_\io), \sg(\ch_\itw))&=&\mcc(\sumjtom \sg(\al_{\io j})F_j, \sumjtom \sg(\al_{\itw j})F_j)=
\\ 
&\numeq{\ref{mtau}} & 
\sumjtom \frac{\dimccj}{\dimcc}\al_{\io\tau(j)}\ovl{\al_{\itw \tau(j)}}=
\\ & \numeq{\ref{alsg}} & 
\sumjtom \frac{\dim(\cc^{\tau(j)})}{\dimcc}\al_{\io\tau(j)}\ovl{\al_{\itw\tau(j)}}=
\\ &=& \sumjtom \frac{\dimccj}{\dimcc}\al_{\io j}\ovl{\al_{\itw j}}=\\
&\numeq{\ref{orthos2t}} & \delta_{\io, \itw}
\end{eqnarray*}
\epf
\subsection{When $\wsgf$ is an algebra map?}
Let $\cc$ be a pivotal fusion category with a commutative character ring $\cfcc$. Recall \cite{scbf8} that $\wcfcc$ is defined as the space of all linear maps $f:\cfcc\ra \comp$. It is a $\comp$-algebra with multiplication defined linearly by
$$
(\mu_\jo\star\mu_\jtw)(\frac{\ch_i}{d_i})=\mu_\jo(\frac{\ch_i}{d_i})\mu_\jtw(\frac{\ch_i}{d_i}).
$$
where $d_i$ is the quantum dimension of the simple objects $X_i$. Since  $\cfcc$ is a commutative ring it follows as above that the set of all algebra homomorphisms $\muj:\cfcc\ra \comp$ is a $\comp$-linear basis for $\cfcc$.

As above, since we assume $\cc^0=\unu_\czcc$ it follows that $\mu_0$ is the unit of $\wcfcc$.
It follows that one can write that
\beq\label{mujmult}
\mujo\star \mujtw=\sumktom\wpp_{k}(\jo, \jtw)\muk
\eeq	
for some scalars $\wpp_{k}(\jo, \jtw)\in \comp$.

From the proof \cite[Theorem 1.1]{scbf8} one has that
\beq\label{wpck}
c^k_{ij}=\frac{\dim(\cc^i)\dim(\cc^j)}{\dim(\cc^k)}{\wdht p}_k(i,j).
\eeq

Note that as above the unit of the $\comp$-algebra $\cfcc$ is given by the quantum dimension function $d:\cfcc\ra \comp,\;\;d(\ch_i)=d_i=\ch_i(1_\cecc)$.

For any $\sg\in \galkq$ recall that we defined a linear map $\wsgf:\wcfcc\ra \wcfcc$ which on the canonical basis given by the algebra homomorphisms $\muj$ is given by $\wsgf(\mu_j):=\sg.\muj=\mu_{\tauj}$.

\bp \label{wsgfalgm}
Let $\cc$ be a  pivotal  fusion category with a commutative Grothendieck ring. 
\bne
\item Then $\wsgf:\widehat{\cfcc}\ra \widehat{\cfcc}$ is an algebra map if and only if
\beq\label{condg}  
\sumktom\sg^{-1}(\wpp_{k}(\jo, \jtw))\al_{ik} =\sumktom\wpp_{k}(\jo, \jtw)\frac{d_i}{\sg^{-1}(d_i)}\al_{ik} 
\eeq
for all $i,k,\jo,\jtw$.
\item 
If $\cc$ is such that $d_i\in \mathbb Z$ for all $i$, then $\wsgf$ is an algebra map if and only if
\beq\label{cond}
\sg^{{-1}}(\wpp_{k}(\jo, \jtw))=\wpp_{k}(\jo, \jtw)
\eeq
for all $k,\jo,\jtw\in \mtc J$.
\item 
It follows that in the case  $d_i\in \mathbb Z$ for all $i$,
$\wsgf$ are algebra maps for all $\sg\in \galkq$ if and only if $\wpp_{k}(\jo, \jtw)$ are rational numbers, for all $\jo, \jtw, k$.
\ene
\ep

\bpf
One has to show that 
$$
\wsgf(\mujo\star \mujtw)=\wsgf(\mujo)\star \wsgf(\mujtw),
$$
for all $\jo, \jtw \in \mtc J$. This is equivalent to 
$$
\wsgf(\mujo\star \mujtw)[\frac{\ch_i}{d_i}]=\wsgf(\mujo)[\frac{\ch_i}{d_i}]\star \wsgf(\mujtw)[\frac{\ch_i}{d_i}],
$$
for any irreducible character $\ch_i$.

As above we have
$\mujo\star \mujtw=\sumktom \wpp_{k}(\jo, \jtw)\muk$ and  
evaluating both sides at $\frac{\ch_{i}}{d_i}$ it follows that
$$
\mujo(\frac{\ch_{i}}{d_i})\mujtw(\frac{\ch_{i}}{d_i})=\sumktom \wpp_{k}(\jo, \jtw)\muk(\frac{\ch_{i}}{d_i})
$$

which can be written as 
\beq\label{5}
\frac{1}{d_i}\al_{i\jo}\al_{i\jtw}=\sumktom \wpp_{k}(\jo, \jtw)\al_{ik}.
\eeq
Applying $\sg^{-1}$ to the above Equation one obtains
\beq\label{6}
\frac{\sg^{-1}(\al_{i\jo}\al_{i\jtw})}{\sg^{-1}(d_i)}=\sumktom\sg^{{-1}}(\wpp_{k}(\jo, \jtw))\sg^{-1}(\al_{ik})
\eeq

\vsk On the other hand, by Equation \eqref{mujmult} note that
\beq\label{7}
\wsgf(\mujo\star \mujtw)=\sumktom \wpp_{k}(\jo, \jtw)\wsgf(\muk)
\eeq
and 
$$
\wsgf(\muj)(\frac{\ch_i}{d_i})=\mu_{\tau(j)}(\frac{\ch_i}{d_i})=\frac{\al_{i\tau(j)}}{d_i}=\frac{\sg(\alij)}{d_i}
$$
Therefore
\begin{eqnarray*}
\wsgf(\mujo\star \mujtw)[\frac{\ch_{i}}{d_i}] 
& = & 
\bigg(\sumktom\wpp_{k}(\jo, \jtw)\wsgf(\muk)\bigg)[\frac{\ch_{i}}{d_i}]=
\\ &=& 
\frac{1}{{d_i}} \bigg(\sumktom \wpp_{k}(\jo, \jtw)\sg(\al_{ik})\bigg)=
\\ &=&  
\sg \bigg( \sumktom \sg^{-1}(\wpp_{k}(\jo, \jtw))\al_{ik}\frac{1}{{\sg^{-1}(d_i)}} \bigg)
\end{eqnarray*}
Moreover,
\begin{eqnarray*}
 \wsgf(\mujo)[\frac{\ch_{i}}{d_i}]  \wsgf(\mujtw)[\frac{\ch_{i}}{d_i}] 
 &=& 
 \frac{1}{d_i^2}\sg(\al_{i\jo})\sg(\al_{i\jtw})=
 \\ &=& 
 \frac{1}{d_i^2} \sg(\al_{i\jo}\al_{i\jtw})=
 \\& 
 \numeq{\ref{5}}&  \frac{1}{d_i^2}\sg\big(\sumktom\wpp_{k}(\jo, \jtw)\al_{ik} d_i\big)=
 \\& = &
 \sg\big(\sumktom\wpp_{k}(\jo, \jtw)\al_{ik}  \frac{d_i}{\sg^{-1}(d_i)^2}
 \big)
\end{eqnarray*}

Therefore $\wsgf$ is an algebra map if and only if
$$  
\sumktom\sg^{-1}(\wpp_{k}(\jo, \jtw))\al_{ik}\frac{1}{\sg^{-1}(d_i)} =\sumktom\wpp_{k}(\jo, \jtw)\al_{ik} \frac{d_i}{\sg^{-1}(d_i)^2}
$$
which can be written
$$  
\sumktom\sg^{-1}(\wpp_{k}(\jo, \jtw))\al_{ik} =\sumktom\wpp_{k}(\jo, \jtw)\frac{d_i}{\sg^{-1}(d_i)}\al_{ik} 
$$
On the other hand if $d_i\in \mz$ then $\sg^{-1}(d_i)=d_i$ and therefore $\wsgf$ is an algebra map if and only if 
$$  
\sumktom\sg^{-1}(\wpp_{k}(\jo, \jtw))\al_{ik} =\sumktom\wpp_{k}(\jo, \jtw)\al_{ik}.
$$
Since $\{\al_{ik}\}$ is an invertible matrix it follows that $\wsgf$ is an algebra map if and only if $\sg^{{-1}}(\wpp_{k}(\jo, \jtw))=\wpp_{k}(\jo, \jtw)$.
\epf

\br
For  weakly integral premodular categories from Equation (4.5) of the proof of \cite[Theorem 1.3]{scbf8} one has 
\beq\label{wprat}
\wdhat{p}_k(i,j)=\sum_{v\in \ca_k}\wdhat{P_v}(s,u)=\sum_{v\in \ca_k}P_v(s,u)=\sum_{v\in \ca_k}\frac{\wdhat{N^v_{su}}\wdhat{d_v}}{\wdhat{d_s}\wdhat{d_u}}.
\eeq
This shows  that in the case of $\cc$ is a weakly integral premodular category  the scalars $\wpp_{k}(\jo, \jtw)$ are rational numbers and therefore $\wsgf$ is an algebra map for all $\sg \in \galkq$.
\er
\bp
Let $\cc$ be an integral premodular category $\cc$. Then $\wsgf$ is an algebra map for all $\sg \in \galkq$.
\ep
\bpf
By the above remark, $\wps_k(i,j)$ are  rational numbers and the result follows from Proposition \ref{wsgfalgm}.
\epf

\subsection{Definition of $\sge$} Let $\cc$ be a pivotal fusion category.
Recall that there are dual bases for the canonical evaluation $$
<, >:\cfcc\ra \cecc \ra \comp, <\ch, z>\sent \ch\circ z
$$
given
\beq\label{db}
\{F_j, \frac{\mtr C_j}{\dim(\cc^j)}\}.
\eeq
By \cite[Theorem 3.12]{scbf8} there is a $\comp$-algebra isomorphism 
$$
\widehat{\cfcc}\xra{\al}\cecc, \;\muj\sent \frac{\mtr C_j}{\dim(\cc^j)}.
$$

\bl \label{eval2}
The isomorphism $\al:\wcfcc \ra \cecc, \muj\sent \frac{\mtr C_j}{\dim(\cc^j)}$ verifies
\beq\label{alev}
<\mu,\ch>=<\ch, \al(\mu)>
\eeq
\el
\bpf 
Indeed, on the canonical  bases one has:
$$
<\muj, \ch_i>=\muj(\ch_i)=\alij=<\ch_i, \frac{\mtr C_j}{\dim(\cc^j)}>.
$$
\epf

Transferring via $\al$ the endomorphism $\wsgf$ on $\cecc$ we obtain an endomorphism 
\beq\label{sgedfn}
\sge:\cecc\ra \cecc,\;\;\frac{\mtr C_j}{\dim(\cc^j)}\sent \frac{\mtr C_{\tau(j)}}{\dim(\cc^{\tau(j)})}.
\eeq
Indeed, $\sge(z)=\al(\wsgf(\al^{-1}(z)))$ for all $z\in \cecc$.
It follows that
$$
\sge(\frac{\mtr C_j}{\dim(\cc^j)})=\al(\wsgf(\al^{-1}(\frac{\mtr C_j}{\dim(\cc^j)})))=\al(\wsgf(\muj))=\al(\mu_{\tau(j)})=\frac{\mtr C_\tauj}{\dim(\cc^\tauj)}
$$
\newcommand{\cctauj}{\cc^{\tau(j)}}
Since $\sge$ is linear it follows that
\beq\label{sgeg}
\sge(\mtr C_j)=\frac{\dimccj}{\dim\cctauj}\mtr C^{\tau(j)}\numeq{\ref{njfus}}\frac{n_{\tau(j)}}{n_j}\mtr C^{\tau(j)}\numeq{\ref{sgnj}}\frac{\sg(n_j)}{n_j}\mtr C^{\tau(j)}
\eeq
for all $j$.

Note that since $\al$ is an algebra isomorphism it follows that  $\wsgf$ is an algebra map if and only if $\sge$ is an algebra endomorphism of $\cecc$.

\bl
If $\cc$ is weakly integral then 
\beq\label{sgeccjwint}
\sge(\mtr C_j)=\mtr C_{\tau(j)}
\eeq
for all $j$.
\el
\bpf
In the weakly integral case one has that $\dim(\ccj)=\dim(\cc^{\tau(j)})$ for all $j$ by Equation \eqref{sgccj2}.
\epf
Note that by Equation \eqref{bnf2} it follows that $c^k_{ij}\in \mbk$ for all $i,j,k$.
\bl 
Let $\cc$ be a pivotal fusion category and $\sg\in \galkq$. Then $\sge$ is an algebra map if and only if
\beq\label{sgegen}
\frac{\sg(n_k)}{n_k}c^k_{\jo,\jtw}=\frac{\sg(n_\jo)}{n_\jo}\frac{\sg(n_\jtw)}{n_\jtw}{c^{\tau(k)}_{\tau(\jo)\tau(\jtw)}}
\eeq
\el
\newcommand{\cctauk}{\cc^{\tau(k)}}
\bpf
Note that $\sge$ is an algebra map if and only if
$$
\sge(\mtr C_\jo \mtr C_\jtw)=\sge(\mtr C_\jo)\sge(\mtr C_\jtw).
$$
On the other hand since $\mtr C_\jo \mtr C_\jtw=\sumktom c^k_{\jo,\jtw}\mtr C_k$ it follows that 
$$
\sge(\mtr C_\jo \mtr C_\jtw)=\sumktom c^k_{\jo,\jtw}\frac{\sg(n_{\tau(k)})}{n_k}\mtr C_{\tau(k)}.
$$
Note also that
\begin{eqnarray*}
\sge(\mtr C_\jo)\sge(\mtr C_\jtw)&=&\frac{\sg(n_\jo)}{n_\jo}\frac{\sg(n_\jtw)}{n_\jtw}\mtr C_{\tau(\jo)}\mtr C_{\tau(\jtw)}=\\ &=&\frac{\sg(n_\jo)}{n_\jo}\frac{\sg(n_\jtw)}{n_\jtw}\bigg(\sumktom c^k_{\tau(\jo),\tau(\jtw)}\mtr C_{k}\bigg)=\\ &=&\frac{\sg(n_\jo)}{n_\jo}\frac{\sg(n_\jtw)}{n_\jtw}\bigg(\sumktom c^{\tau(k)}_{\tau(\jo),\tau(\jtw)}\mtr C_{\tau(k)}\bigg).
\end{eqnarray*}
Comparing the above two Equations we obtain that $\sge$ is an algebra map if and only if $\frac{\sg(n_k)}{n_k}c^k_{\jo,\jtw}=\frac{\sg(n_\jo)}{n_\jo}\frac{\sg(n_\jtw)}{n_\jtw}{c^{\tau(k)}_{\tau(\jo)\tau(\jtw)}}$.
\epf

\bc\label{ccwint}
If $\cc$ is an integral fusion category then $\sge$ is an algebra map if and only if $\sg(c^k_{\jo,\jtw})=c^k_{\jo,\jtw}$ for all $\jo,\jtw,k$.
\ec
\bpf
If $\cc$ is integral then $\sg(n_j)=n_j$ and the result follows from Proposition \ref{strcint}.
\epf
\bl\label{eval} 
Let $\cc$ be a pivotal fusion category. For all $\ch \in \cfcc$ and $z\in \cecc$ one has
\beq\label{evchz}
<\ch, \sge(z)>=<\sgf(\ch), z>.
\eeq
\el
\bpf 
Since $<,\;>$ is bilinear it is enough to verify the above identity on the basis $\{\ch_i\}$ of $\cfcc$ and $\{\frac{\mtr C_j}{\dim(\cc^j)}\}$ of $\cecc$. Indeed, one has
$$
<\ch_i, \sge(\frac{\mtr C_j}{\dimccj})>=<\ch_i, \frac{\mtr C_\tauj}{\dimcctauj}>=\mu_\tauj(\ch_i)=\al_{i\tau(j)}.
$$
On the other hand,
$$
<\sg.\ch_i, \frac{\mtr C_j}{\dimccj}>=<\sum_{l=0}^m\sg(\al_{il})F_l,\;\frac{\mtr C_j}{\dimccj}>\numeq{\ref{db}}
\al_{i\tau(j)}.
$$
\epf
\bl \label{fchi} 
Let $\cc$ be a pivotal fusion category. For any $\ch_i \in \cfcc$ an irreducible character and all $z, z'\in \cecc$ one has
\beq\label{zz}
d_i\ch_i(zz')=\ch_i(z)\ch_i(z').
\eeq
\el
\bpf 
If $z=\sumltom z_lE_l$ and $z'=\sumltom z'_lE_l$ then $zz'=\sumltom z_lz'_lE_l$
Then $\ch_i(z)\ch_i(z')=z_iz'_id_i^2$ and $d(\ch_i)\ch_i(zz')=d_i^2z_iz'_i=\ch_i(z)\ch_i(z')$.
\epf

\bl\label{reverse} 
Let $\cc$ be a pivotal fusion category. Let $\ch \in \cfcc$ and $d \in \comp\setminus\{0\}$. Then 
$$
\ch(zz')=\frac{1}{d}\ch(z)\ch(z')
$$
for all $z, z'\in \cecc$ if and only if
$$
\ch=\frac{d}{d_i}\ch_i
$$
for some irreducible character $\ch_i$ of $\cc$.
\el
\bpf
Suppose that $\ch=\sum_{i \in \cca}\al_i\ch_i$ where $\cca\subseteq \{0,1,\dots ,m\}$ is the set of all irreducible characters that have a non-zero coefficient $\al_i\neq 0$.

Let $z=E_\io$ and $z'=E_\itw$. It follows that $0=\frac{\al_\io d_\io \al_\itw d_\itw}{d}$ which is a contradiction. Therefore $\cca$ is a single element and without loss of generality we may suppose $\cca=\{i\}$. Therefore $\ch=\al \ch_i$ and for $z=z'=1_{\cecc}$ this gives that $\al=\frac{d}{d_i}$.
\epf
\bp \label{sgebyev}
Let $\cc$ be a pivotal fusion category and $\sg \in \galkq$. Then 
$\sge$ is an algebra map if and only if there is a permutation $\eta$ of the indices  $\{0,\dots ,m\}$ such that 
\beq\label{rev2}
\sgf(\ch_i)=\frac{d_i}{d_{\eta(i)}}\ch_{\eta(i)}.
\eeq 
for all $i\in  \{0,\dots ,m\}$.
\ep
\bpf 
One has that $\sge$ is multiplicative if and only if 
$$\sge(zz')=\sge(z)\sge(z')$$ for all $z,z'\in \cecc$. Since the evaluation form is non-degenerate it follows that this is equivalent to
$$
<\ch_i, \sge(zz')>=<\ch_i, \sge(z) \sge(z')>,
$$
for all $z,z'\in \cecc$ and all $i$.

Note that
\begin{eqnarray*}
<\ch_i, \sge(zz')>&\numeq{\ref{evchz}}&<\sgf(\ch_i), zz'>
\end{eqnarray*}
and
\begin{eqnarray*}
<\ch_i, \sge(z)\sge(z')>&\numeq{\ref{zz}}&\frac{1}{d_i}<\ch_i, \sge(z)><\ch_i, \sge(z')>=\\&\numeq{\ref{evchz}}&\frac{1}{d_i}<\sgf(\ch_i),z><\sgf(\ch_i), z'>
\end{eqnarray*}
Therefore $\sge$ is an algebra map if and only if 
\beq\label{sge2}
<\sgf(\ch_i), zz'>\numeq{\ref{evchz}}\frac{1}{d_i}<\sgf(\ch_i),z><\sgf(\ch_i), z'>
\eeq
for all $\ch_i\in \irr(\cc)$ and all $z,z'\in \cecc$.

"$\implies$"  Suppose that $\sge$ is an algebra map. Equation \eqref{sge2} implies by Lemma \ref{reverse} that
that $\sgf(\ch_i)$ is a scalar of an irreducible character of $\cc$. Therefore for all $i\in \{0,\dots ,m\}$ there is an index $\eta_i$ in $\{0,\dots ,m\}$ such that $\sgf(\ch_i)=\frac{d_i}{d_{\eta(i)}}\ch_{\eta(i)}$. Since $\sgf$ is bijective it follows that $\eta$ is a permutation of $ \{0,\dots ,m\}$.

"$\impliedby$" Conversely, if there is a permutation $\eta$ of the indices  $\{0,\dots ,m\}$ such that $\sgf(\ch_i)=\frac{d_i}{d_{\eta(i)}}\ch_{\eta_i}$ then by Lemma \ref{fchi} it follows that the Equation \eqref{sge2} is satisfied. Thus $\sge$ is an algebra map in this case.
\epf

\bp\label{irrsg} 
Let $\cc$ be a pivotal fusion category and $\sg\in \galkq$. If $\sge$ is an algebra map then $\sge(E_i)=E_{\eta^{-1}(i)}$ for all $i$.
\ep
\bpf
Since  $\sge$ is an algebra isomorphism it permutes the central primitive $E_i$. Thus one can write $ \sge(E_i)=E_{\ro(i)}$ for some permutation $\ro=\ro_{\sg}$ of $\mtc I$.

Recall that by Equation \eqref{evchz} oen has
$$
<\ch, \sge(z)>=<\sgf(\ch), z>,
$$
for all $\ch \in \cfcc$ and $z\in \cecc$.

For $\ch=\ch_i$ and $z=E_s$ in the above equation it follows that
$$
<\ch_i, \sge(E_s)>=<\ch_i, E_{\ro(s)}>=\delta_{i,\ro(s)}d_i.
$$ 
On the other hand, since $\sge$ is an algebra map it follows by Equation \eqref{rev2} that 
$$
<\sgf(\ch_i), E_s>=<\frac{d_i}{d_{\eta(i)}}\ch_{\eta(i)}, E_s>=\frac{d_i}{d_{\eta(i)}}<\ch_{\eta(i)}, E_s>=\delta_{\eta(i),s}d_i.
$$

This shows that $\eta^{-1}(i)=\ro(i)$ for all $i$, i.e. $\eta^{-1}=\ro$.
\epf
\bp
Let $\cc$ be a pivotal fusion category. Suppose that $\sge$ is an algebra map for some $\sg\in \galkq$. With the above notations it follows that 
\beq\label{etaist}
\eta(i^*)=\eta(i)^*
\eeq
and 
\beq\label{diist1}
d_id_{i^*}=d_{\eta(i)}d_{\eta(i)^*}
\eeq
for all $i$.
\ep
\bpf
Applying $\sgf\ot \sgf$ to Equation \eqref{clseqcc} one obtains by Equations \eqref{sgnj} and \eqref{sgfj} that
\beq\label{sgclseq}
\sumjtom n_jF_j\ot F_j=\sumitom \sgf(\ch_i)\ot \sgf(\ch_i^*).
\eeq
Since $\sge$ is an algebra map it follows that Equation \eqref{sgclseq} can be written as
\beq\label{sgealg}
\sumjtom n_jF_j\ot F_j=\sumitom \frac{d_i}{d_{\eta(i)}}\ch_{\eta(i)}\ot \frac{d_{i^*}}{d_{\eta(i^*)}}\ch_{\eta(i^*)}.
\eeq
Using again  Equation \eqref{clseqcc} it follows that 
$
\eta(i^*)=\eta(i)^*
$
and 
$
d_id_{i^*}=d_{\eta(i)}d_{\eta(i)^*}
$
for all $i$.
\epf

\br\label{w-intsge}
In the weakly integral case the relation ${d_i}^2={d_{\eta_i}}^2$ follows also from Equation \eqref{selfm}.
\er

\bc\label{psu}
If $\cc$ is a spherical fusion category $\sg_E$ is an algebra map then with the above notations one has $d_i=\pm d_{\eta(i)}$ for all $i$.
If $\cc$ is  a pseudo-unitary fusion category and $\sg_E$ is an algebra map then with the above notations one has $d_i=d_{\eta(i)}$ for all $i$.
\ec

\bpf If $\cc$ is spherical then $d_i=d_{i^*}$ are real numbers and therefore $d_i^2=d_{\eta(i)}^2$. Moreover, if $\cc$ is  a pseudo-unitary fusion category then $d_i>0$ in this case and therefore $d_i^2=d_{\eta(i)}^2$ implies $d_i=d_{\eta(i)}$.
\epf

\bc\label{sgdi}
Let $\cc$ be a pivotal fusion category and $\sg\in \galkq$ such that $\sge$ is an algebra map. Then  
$$
\sg(d_i)=d_i
$$ 
for all $i\in  \{0,\dots ,m\}$.
\ec
\bpf 
Note that by applying dimensions to Equation \eqref{rev2} it follows that $d(\sgf(\ch_i))=d_i$ for all $i\in  \{0,\dots ,m\}$. On the other hand, $d(\sgf(\ch_i))=\sg(\al_{i0})=\sg(d_i)$ for all $i\in  \{0,\dots ,m\}$.
\epf

\bc \label{qint}
Let $\cc$ be a pivotal fusion category and suppose that $\sge$ is an algebra map for all $\sg$ and that the extension $\mbq \subseteq \mbk$ is Galois. In this case $d_i\in \mz$  for all $i$.
\ec
\bpf
$d_i$ are rational from previous Corollary and algebraic integers therefore they are integers.
\epf
\bp 
Let $\cc$ be a pivotal fusion category and $\sg\in \galkq$. If $\sge$ is an algebra map then $\sgf(\lam)=\lam$, i.e. $\tau(0)=0$.
\ep
\bpf 
By \cite[Equation (6.8)]{scalg} one has
$$
F_0=\lam=\frac{1}{\dimcc}\big(\sumitom d_{i^*}\ch_i \big)
$$
and therefore
\beq\label{sglam2}
\sgf(\lam)=\frac{1}{\dimcc}\big(\sumitom d_{i^*}\sgf(\ch_i)\big)
\eeq

If $\sge$ is an algebra map then there is a permutation $\eta=\eta_\sg$ such that 
$ \sgf(\ch_i)=\frac{d_i}{d_{\eta(i)}}\ch_{\eta(i)}$ for all $i$.
Then Equation \eqref{sglam2} implies that
\beq\label{sglam3}
\sgf(\lam)=\frac{1}{\dimcc}\big(\sumitom d_{i^*} \frac{d_i}{d_{\eta(i)}}\ch_{\eta(i)}\big)\numeq{\ref{diist1}}\frac{1}{\dimcc}
\big(
\sumitom \frac{d_{\eta(i)}d_{\eta(i)^*}}{d_{\eta(i)}}\ch_{\eta(i)}
\big)=\lam
\eeq
On the other hand by our assumption we have $\lam=F_0$ and therefore $\sgf(\lam)=F_{\tau^{-1}(0)}$.
\epf

\subsection{Proof of Theorem \ref{int-main2}}
\bpf
Note that $\dimccj$  are integers since $\cc$ is integral. The proof of the theorem follows from Proposition \ref{wsgfalgm} by noticing that $\widehat{p}_k(\jo,\jtw)$ are rational numbers by Equation \eqref{wpck}.
\epf
\bp 
Let $\cc$ be a pivotal fusion category and $\sg \in \galkq$. If $\sge$ is an algebra map then one has that 
\beq\label{g2}
\alij\ovl{\bigg(\frac{d_i}{d_{\eta(i)}}\bigg)}=\ovl{\bigg(\frac{n_j}{\sg^{-1}(n_j)}\bigg)}\sg^{-1}(\al_{\eta(i) j}).
\eeq
for all $i, j$.
\ep
\bpf
By Equation \eqref{fj} one has that
$$
F_j=\frac{1}{n_{j}}\big(\sumitom \ovl{\muj(\ch_{i})}\ch_{i}\big)
$$ 
Applying $\sgf$ to this formula one has that 
$$
\sgf(F_j)=\frac{1}{n_{j}}\big(\sumitom \ovl{\muj(\ch_{i})}\sgf(\ch_{i})\big)\numeq{\ref{rev2}}\frac{1}{n_{j}}\big(\sumitom \ovl{\muj(\ch_{i})}\frac{d_i}{d_{\eta(i)}}\ch_{\eta(i)}\big).
$$ 
On the other hand by Proposition \ref{sgfonfj} one has
$$
\sgf(F_j)=F_{\tau^{-1}(j)}=\frac{1}{n_{\tau^{-1}(j)}}\big(\sumitom \ovl{\mu_{\tau^{-1}(j)}(\ch_{i})}\ch_{i}\big)=\frac{1}{n_{\tau^{-1}(j)}}\big(\sumitom \ovl{\mu_{\tau^{-1}(j)}(\ch_{\eta(i)})}\ch_{\eta(i)}\big)
$$
Since $n_j=\sg(n_{\tau^{-1}(j)})$ by Equation \eqref{sgnj}, from the two formula for $\sgf(F_j)$ one obtains that:
$$
\ovl{\muj(\ch_{i})}\frac{d_i}{d_{\eta(i)}}=\frac{n_j}{\sg^{-1}(n_j)}\ovl{\mu_{\tau^{-1}(j)}(\ch_{\eta(i)})}
$$
for all $i, j$. This can also be written as:
$$
\alij\frac{\ovl{d_i}}{\ovl{d_\eta(i)}}=\ovl{\frac{n_j}{\sg^{-1}(n_j)}}\al_{\eta(i) \tau^{-1}(j)}=\ovl{\frac{n_j}{\sg^{-1}(n_j)}}\sg^{-1}(\al_{\eta(i) j})
$$
for all $i, j$.
\epf
Note that Equation \eqref{g2} can also be written as 
$$
\sg(\alij)=\al_{i\tau(j)}=\ovl{\bigg(\frac{n_j}{\sg^{-1}(n_j)}\bigg)}\ovl{\bigg(\frac{d_i}{d_{\eta(i)}}\bigg)}\al_{\eta(i)j}.
$$
which in the weakly-integral case can be written as 
$$
\sg(\alij)=\al_{i\tau(j)}=\al_{\eta(i)j}
$$
for all $i,j$.
\section{Zeros in the character table}\label{zeros}
A classical result of Burnside in character theory states that for any  
irreducible character $\ch$ of a finite group $G$ with $\ch(1)>1$ there is some
$g \in G$ such that $\ch(g)=0$, see \cite[Chapter 21]{bz}. This result was generalized in \cite{gnn} to weakly integral modular categories. Recall \cite{eno-adv} that a fusion category $\cc$ is called {\it weakly integral} if its Frobenius-Perron dimension is an integer. In this case the Frobenius-Perron dimension of every simple object of $\cc$ is the square root of an integer \cite{eno-annals}.

In this section we extend further the above result to  integral  fusion categories with rational structure constants $c^k_{ij}$. The proof of this result goes along the same lines as the classical result proven in \cite[Chapter 21]{bz} with a slight modification concerning some AM-GM inequality. For the sake of completeness we include a sketch the proof of the result below.

We keep the same notations as in the previous sections.
\bpf 
For any $0\leq i \leq m$ denote $\mtc T_i:=T(\ch_i)=\{j\;|\; \alij= 0\}$ and $\cd_i:=\mtc J\setminus (\mtc T_i\cup \{0\})$.
\vsk
The first orthogonality relation from Equation \eqref{611} can be written as:
 \beq\label{orthofes}
\sum_{ j=0}^m\frac{|\ch_i({\mtr C_j})|^2}{\dim(\cc^j)}=\dimcc.
\eeq

Since for $j=0$ one has $\mtr C_0=\unue$ this can be written as
$$
\dimcc=d_i^2+\sum_{j \in \cd_i}\frac{|\ch_i({\mtr C_j})|^2}{\dim(\cc^j)}
$$
which gives that 
\beq\label{gnn1}
1=\frac{\dimcc}{d_i^2}-\sum_{j \in \cd_i}\frac{|\ch_i({\mtr C_j})|^2}{d_i^2\dim(\cc^j)}
\eeq
On the other hand note that
$$
\dimcc=\sumjtom \dim(\cc^j)=1+\sum_{j\in \mtc T_i}\dim(\cc^j)+\sum_{j\in \cd_i}\dim(\cc^j)
$$
Therefore Equation \eqref{gnn1} can be written as:
\beq\label{gnn2}
1=\frac{1+\sum_{j\in \mtc T_i}\dim(\cc^j)}{d_i^2}-(\sum_{j \in \cd_i}\frac{|\ch_i({\mtr C_j})|^2}{d_i^2\dim(\cc^j)}-\sum_{j\in \cd_i}\frac{\dim(\cc^j)}{d_i^2}
)
\eeq
Thus in order to finish the proof it is enough to show that
\beq\label{gnn2}
\sum_{j \in \cd_i}\frac{|\ch_i({\mtr C_j})|^2}{d_i^2\dim(\cc^j)}-\sum_{j\in \cd_i}\frac{\dim(\cc^j)}{d_i^2}\geq 0,
\eeq
since then it follows that $\frac{1+\sum_{j\in \mtc T_i}\dim(\cc^j)}{d_i^2}\geq 1$, i.e. $1+\sum_{j\in \mtc T_i}\dim(\cc^j)\geq d_i^2$. Since $d_i>1$ it follows that $\mtc T_i\neq\emptyset$.
\vsk 
The inequality from Equation\eqref{gnn2} can be written as
\beq\label{gnn3}
\frac{1}{\sum_{j \in \cd_i}\dim(\cc^j)}(\sum_{j \in \cd_i}\frac{|\ch_i({\mtr C_j})|^2}{\dim(\cc^j)})\geq 1.
\eeq
\noindent On the other hand the weighted AM-GM inequality gives that 
\beq\label{gnn4}
\frac{1}{\sum_{j \in \cd_i}\dim(\cc^j)}\bigg(\sum_{j \in \cd_i}\frac{|\ch_i({\mtr C_j})|^2}{\dim(\cc^j)}\bigg)\geq (\prod_{j\in \cd_i} \bigg( \big(\frac{|\ch_i({\mtr C_j})|^2}{\dim(\cc^j)^2}\big)^{\dim(\cc^j)} \bigg)^{\frac{1}{|\cd_i|}},
\eeq
where $|\cd_i|:=\sum_{j\in \cd_i}\dimccj$.
Note that Equation \eqref{alsg} implies that the set $\cd_i$ is stable under the Galois group $\galkq$. This in turn implies that the product  
$$P_i:=\prod_{j\in \cd_i}(\frac{|\ch_i({\mtr C_j})|^2}{\dim(\cc^j)^2})^{\dim(\cc^j)}$$ is fixed by the Galois group $\galkq$ 
 since $\dim(\cc^{\tau(j)})=\dim(\cc^j)$ by Equation \eqref{sgccj2}.
It follows that $P_i$ is a rational number. On the other hand each factor of this product is an algebraic integer (since $\dimccj \in \mathbb Z_{_{>0}}$) and therefore the entire product is an integer. Since it is positive it follows it is greater or equal to $1$.
\epf
In analogy with group representations, we call a fusion category {\it perfect} if it has no non-trivial invertible objects. Next result generalizes a well-known result of Brauer from group representation theory to integral perfect braided fusion categories.
\bt
Let $\cc$ be a weakly integral braided fusion category. Then $\cc$ is a perfect fusion category if and only if the following identity holds:
\beq
\sumjtom \mtr C_j=\frac{\dimcc}{\prod_{j=0}^m\dim(\cc^j)}(\prod_{j=0}^m \mtr C_j)
\eeq
\et
\bpf
If $\cc$ is perfect it is enough to show that 
$$
\ch_i(\sumjtom \mtr C_j)=\frac{\dimcc}{\prod_{j=0}^m\dim(\cc^j)}\ch_i(\prod_{j=0}^m \mtr C_j)
$$
for any irreducible character $\ch_i \in \irr(\cc)$. For $\ch_i=\ch_0=\epsu$ one obtains equality by dimension argument. For $\ch_i\neq \ch_0$ both terms above are zero. Indeed note that  $\lam_\cc=\frac{1}{\dimcc}(\sumjtom \mtr{C}_j)$ and therefore $\ch_i(\sumjtom \mtr C_j)=0$. On the other hand the right hand side is zero by the vanishing theorem since $\ch_i$ is the character of a non-invertible object. Conversely, note that if $M_i$ is an invertible object it follows that $\ch_i(\mtr C_j)\neq 0$ for any character $j$. Indeed one can write $\ch_i=\sumjtom \al_jF_j.$ Since there is $n\geq 0$. such that $\psi^n=\epsu=\sumjtom F_j$ it follows that $\al_j\neq 0$ for all $j$. On the other hand note that, by Equation \eqref{normcj} one has 
$\al_j=\frac{1}{\dim(\cc^j)}\ch_i(\mtr{C}_j)$.
\epf
\subsection{The modular case} Recall that the $S$-matrix of a braided pivotal fusion category is defined as $S_{XY}:=Tr(C_{X,Y}C_{Y,X})$. Then if $S_{ij}:=S_{M_i, M_j}$ it follows $S_{ij}=S_{ji},\;\; S_{i^*j}=\ovl{S_{ij}}=S_{ij^*}$
for all $i,j$.

As usually, by $\mathbb Q(S)$ is denoted  the field obtained by adjoining the $S$-matrix entries $S_{ij}$ to $\mathbb Q$. By \cite[Theorem 8.14.7]{EGNO15} one has that  $\mathbb Q(S)$ is also contained in a cyclotomic extension.

By \cite[Theorem 6.2]{ccc} in the modular case one has that
$\dim(\cc^k)=d_k^2$ and $\al_{ij}=\frac{S_{ij}}{d_j}$. The first orthogonality from Equation \eqref{611} can be written in this case as
\beq\label{611p}
\sum_{ k=0}^m S_{ik} S_{m^*k}=\delta_{i,m}\dimcc.
\eeq
Since $S_{0j}=d_j$ one has that
$$
\mathbb Q\subseteq \mathbb K=\mathbb Q(\frac{S_{ij}}{d_j})\subseteq \mathbb Q(S)=\mathbb Q(S_{ij}).
$$
By the Fundamental Theorem of Galois theory, since the extension $\mathbb Q\subseteq \mbk_\cc=\mathbb Q(\frac{S_{ij}}{d_j})$ is Galois it follows that there is a group epimorphism 
$
G\xra{\pi} G_\cc.
$
From here it follows that the two actions on $\mu_j$'s are compatible with the epimorphism $\pi$.
\bibliographystyle{alpha}
\bibliography{bitfc17}
\ed
\newpage
\br\label{sgeunitar} The condition $\tau(0)=0$ can also be verified as following. If $\sge$ is multiplicative and bijective it also follows that $\sge$ is unitary.
Indeed, one has $\mtr C_0=\unue$ and therefore by Equation \eqref{sgedfn} one has  $\sge(\mtr C_0)=\frac{\mtr C_{\tau(0)}}{\dim(\cc^{\tau(0)})}$. Thus $\sge(\unue)=\unue \iff \tau(0)=0.$
\er
\blue{Since by our choice one has $\mu_0=d$ it follows that  that $d(\sgf(\ch_i))=\sg(\al_{i0})=\sg(d_i)$ for all $i$.
\\
Then 
$$
(\sg.\mu_0)(\ch_i)=\sg(\mu_0(\ch_i))=\sg(d_i).
$$
}
\newcommand{\wtf}{\widehat{F}}
\newpage
\subsubsection{$\sgf(\lam)$ in two ways for the integral case}

\bl 
Let $\cc$ be a pivotal fusion category. With the above notations one has
\beq\label{sglam}
\sgf(\lam)=\frac{1}{\dimcc}\big(\sumitom \ovl{\sg^{-1}(d_i)}\ch_{i}\big).
\eeq
\el
\bpf
Note that since $\ch_i\lam=d_i\lam$ it follows that $\sgf(\ch_i)\sgf(\lam)=d_i\sgf(\ch_i)$ for all $i$.

It follows that
$$\ch\sgf(\lam)=\sg^{-1}(d(\ch))\sgf(\lam)$$ 
\blue{$\sg^{-1}(d(\ch))$ does not make sense if $\ch$ is  any complex char, values out of $\mbk$. not alg coeff!}

\blue{Need to show that $\mbk$ is a splitting field for $\cfcc$.}
which shows that $\sgf(\lam)$ is the primitive central idempotent corresponding to the morphism $\sg^{-1}\circ \mu_0:\cfcc\ra \comp$.

From the formula \ref{fj} one has $F_j:=\frac{1}{n_{j}}\big(\sumitom \ovl{\muj(\ch_{i})}\ch_{i}\big)$ which gives that
\begin{eqnarray*}
F_{\tau^{-1}(0)}&=&\frac{1}{n_{\tau^{-1}(0)}}\big(\sumitom \ovl{(\sg^{-1}\circ d)(\ch_{i})}\ch_{i}\big)=
\\&=&
\frac{1}{n_0}\big(\sumitom \ovl{\sg^{-1}(d_i)}\ch_{i}\big)=
\\&=&
\frac{1}{\dimcc}\big(\sumitom \ovl{\sg^{-1}(d_i)}\ch_{i}\big)
\end{eqnarray*}
Here we have used that $n_0=n_{\tau^{-1}(0)}=\dimcc$  
\blue{w-int needed!}.
\epf

\newpage
\blue{this remark is wrong, $\al_i$ may be negative!}
\br
Therefore $\delta_{i,\etas}d_i=\mcc(\sg^{-1}(\ch_i), \ch_s)d_s$ for all $s$ which shows that $\sg^{-1}(\ch_i)=\al_i\ch_{\eta^{-1}(i)}$ is a scalar multiple of an irreducible character. 
On the other hand 
\begin{eqnarray*}
\al_i^2&=&\mcc(\sg^{-1}(\ch_i), \sg^{-1}(\ch_i))=\mcc(\sumjtom \sg^{-1}(\al_{ij})F_j, \sumjtom \sg^{-1}(\al_{ij})F_j)
\\ &\numeq{\ref{mtau2}} & \sumjtom \dimccj\sg^{-1}(\alij)\ovl{\sg^{-1}(\alij)}=\sg^{-1}\big(\sg(\dimccj)\sumjtom \alij\al_{i^*j}\big) \numeq{\ref{orthos2t}} 
\end{eqnarray*}
 One has from  above $d_i=\al_i d_{\eta^{-1}(i)}$ therefore $\al_i>0$, since $d_i>0$ for any $i$. This  implies $\al_i=1$.
\er
\newpage
\br
Suppose that $\sge$ is an algebra automorphism. It follows that $(\sg^{-1})_E$ is also an algebra map and therefore there is a permutation $\eta_{\sg^{-1}}$ such that
$$
(\sg^{-1})_F(\ch_i)=\frac{d_i}{d_{\eta_{\sg^{-1}}(i)}}\ch_{\eta_\sg^{-1}(i)}
$$
for all $i$. Since $\sgf\circ (\sg^{-1})_F=\id$ it follows that $\eta_{\sg^{-1}}=\eta_{\sg}^{-1}$.
\er
\blue{This remark not used anywhere!}
\bl \label{eval1}
One has 
\beq\label{evsgs}
<\wsgf(\mu),\ch>=<\mu, \sgf(\ch)>.
\eeq 
\el
\bpf 
Indeed, on the canonical  bases one has:
$$
<\wsgf(\muj), \ch_i>=\mu_{\tau(j)}(\ch_i)=\al_{i\tau(j)}\numeq{\ref{alsg}}\sg(\alij)=<\muj, \sum_{l=0}^m\sg(\al_{il})F_l>=<\muj, \sgf(\ch_i)>.
$$
\epf
\section{Connections with the center}
\subsection{Relation with the center in general}
Let $\cc$ be a pivotal fusion category with a commutative Grothendieck ring $\cfcc\simeq \mtr{Gr}_{\kk}(\cc)$. Let also $F_0, F_1, \dots ,F_m$ be the primitive central idempotents of $\cfcc$ and  $\mu_0, \mu_1, \dots \mu_m$ be their corresponding characters on $\cfcc$. Therefore
$$\mui:\cfcc\ra \comp,\;\mui(F_j)=\delta_{i,j}.$$
We also denote by $\cc^0, \cc^1,\dots \cc^m$ the conjugacy classes of $\cc$ corresponding in this order to the primitive idempotents $F_0, F_1, \dots F_m$.

Let also $M_0, M_1, \dots M_m$ be a complete set of representatives for the isomorphism classes of simple objects of $\cc$. As above, without loss of generality we may assume that $M_0=\unu$ is the unit object of $\cc$.

As above we denote by $V_0, V_1,\dots V_r$  a complete set of  simple objects of $\czcc$ and by $\we_0, \we_1, \dots \we_r$ their associated primitive idempotents in $\mtr{CE}(\czcc)$. We denote also by $\wch_0, \wch_1, \dots \wch_r$ the characters associated to $V_0, V_1, \dots V_r$ and let $\wzd_0,\;\wzd_1,\;\dots \wzd_r$ be their quantum dimensions. Therefore $\wzd_s=\wch_s(\unu)$ for all $s$.
We may also assume that $V_0=\unu$ is the unit object of $\czcc$. \blue{Without loss of generality we may also suppose that $V_i=\cc^i$ for any $0\leq i \leq m$.} Since the Drinfeld map $F_Q:\cf(\czcc)\ra \mtr{CE}(\czcc)$ is bijective it follows that $\wf_j:=F_Q^{-1}(\we_j)$ is a complete set of primitive orthogonal idempotents of $\cf(\czcc)$.

Let also $F:\czcc\ra \cc$ be the forgetful functor. It is well known, see \cite[Lemma 8.49]{eno-annals}, that the induced map $\res=F_*:\cf(\czcc)\ra \cfcc$ is surjective.

Moreover, Ostrik showed in \cite[Theorem 2.13]{O3} that for any primitive idempotent $F_j\in \cfcc$ of a spherical category $\cc$ there is a unique primitive idempotent $\wtf_{\sgj}\in \cf(\czcc)$ whose restriction to $\cfcc$ is $F_j$ and moreover $V_{\sgj}=\cc_{\sg(j)}$ is a conjugacy class of $\cc$ with 
$$
\dim(V_{\sgj}) =\dim(\ccj),\;\text{i.e}\; \wzd_{\sg(j)}=\wzd_j.
$$
Thus 
$F_*(\wf_{\sg(j)})=F_j,\;\tetx{for all}\; 0\leq j \leq m.$ Note also that $F_*(\wf_s)=0$ for $s\neq \sg(j)$ for some $j$. We denote by $\sz_{ij}$ the $S$-matrix of the braided category $\czcc$.

\subsection{Restriction coefficients}
One can write that
\beq\label{rti}
\wch_{\tilde i}=\sum_{s\in \mtc R_{\tilde i}}R^{\tilde i}_s\ch_s,\;\;\text{with} \in \mz_{\geq 0}.
\eeq
It follows from here that 
$$
d_{\tilde i}=\sum_{s\in \mtc R_{\tilde i}}R^{\tilde i}_sd_s\in \mbq(d_s)\subseteq \mbk.
$$
Note that $\al_{i0}=d_i\in \mbk$.
\subsection{Relation between fields and Galois Groups}

Let $\tilde G:=\gal(\mathbb Q(S)/\mbq)$.
\bl For any pivotal fusion category one has that 
$$\mbqs\subseteq \mbk.$$
\blue{Just one part of $\sij$ are in $\mbk$ the lemma might be incorrect or the proof not complete.}
\el
\bpf
As above we may suppose that 
$F_*(\wf_{\sg(j)})=F_j,\;\tetx{for all}\; 0\leq j \leq m.$ Note also that $F_*(\wf_s)=0$ for $s\neq \sg(j)$ for some $j$. It follows that
$$\muj\circ F_*=\wmuj, \;\tetx{for all}\; 0\leq j \leq m.$$
Following \cite[Example 6.14]{scalg}, inside $\cf(\czcc)$ one can write that
$$
\wch_{\tilde i}=\sum_{l=0}^{r} \frac{S_{\tilde i l}}{\wzd_l} \wf_l, 
$$
for all $0\leq i\leq m$.

Applying the morphism $F_*:\cf(\czcc)\ra \cfcc$ induced by the forgetful functor $F$ one obtains that
\beq\label{resti}
F_*(\wch_{\tilde i})=\sum_{j=0}^{r} \frac{S_{\tilde i j}}{\wzd_{j}}F_*(\wf_j)=\sum_{j=0}^m\frac{S_{\tilde i \sg(j)}}{\widehat{d}_{\sg(j)}}F_j
.
\eeq

On the other hand by Equation \eqref{rti} it follows that
\beq\label{rti2}
F_*(\wch_{\tilde i})=\sum_{s\in \mtc R_{\tilde i}}R^{\tilde i}_s\ch_s=\sum_{s\in \mtc R_{\tilde i}}R^{\tilde i}_s\bigg(\sumjtom \al_{sj}F_j\bigg)=\sumjtom \bigg(\sum_{s\in \mtc R_{\tilde i}}R^{\tilde i}_s\al_{sj}\bigg)F_j.
\eeq

Comparing the two Equations it follows that
\beq\label{rtic}
\frac{S_{\tilde i j}}{\wzd_{j}}=\sum_{s\in \mtc R_{\tilde i}}R^{\tilde i}_s\al_{sj},
\eeq
for all $\tilde i$. This shows that $S_{\tilde i j}\in \mbk$.

This implies that $\muj(\ch_t)=\frac{\sz_{s \sg(j)}}{d_{\sg(j)}}=\wmu_{\sgj}(\wch_{\tilde t}),$ thus $\wmu_{\sgj}\in \ca_j$.
\epf
\blue{Read also the Corollary 8.53 from eno-annals again. understand it completely!}
\subsection{} For the rest of this section we suppose that $\cc$ is a braided spherical fusion category, i.e a premodular category. Since $\czcc$ is a modular tensor category in this case it follows by \cite[Theorem 3.1.12]{BaKi} that the irreducible characters of $\cf(\czcc)$ are indexed by the simple object of $\czcc$. More precisely, if $V_i$ is a simple object of $\czcc$ then
$$
\wmu_i:\cf(\czcc)\ra \comp,\; \wmui([V_j])=\frac{\sz_{ij}}{d_i}
$$
is an algebra homomorphism.


By \cite[Section 2.10]{dgno2} there is also a braided tensor functor 
\beq\label{iota}
\iota:\cc \hookrightarrow \czcc, X\sent (X, c_{X,-}).
\eeq
that is fully faithful.
It follows that $\iota(M_i)\simeq V_{\tilde{i}}$ for  some $0\leq \tilde{i}\leq r$.  Note that $\{\tilde 0, \tilde 1, \dots ,\tilde m\}\cap \{0,1,\dots,  m\}=\{0\}$. Indeed, for $i>1$ $\iota(M_i)$ cannot be a conjugacy class since $M_i=F(\iota(M_i))$ does not contain the unit object $\unu$ of $\cc$. Since $F\circ \iota=\id_\cc$ note that 
$F_*(\wch_{\tilde t})=\ch_t,\;\text{for all}\; 0 \leq t \leq m$.

Consider the inclusion $\cfcc\subseteq \cf(\czcc)$ induced by the inclusion functor of \eqref{iota}. For any $0\leq j \leq m$ we define a class of characters
$$
\ca_j:=\{\wmui\in \widehat{\cf(\czcc)}\;|\;\wmui|_{\cfc}=\muj\}.
$$
\bp\label{sigma} Let $\cc$ be a premodular category.
With the above notations one has $\sg(j)\in \ca_j$ for any  $0\leq j \leq m$. 
\ep
\bpf 
Following \cite[Example 6.14]{scalg}, inside $\cf(\czcc)$ one can write that
$$
\wch_{\tilde i}=\sum_{l=0}^{r} \frac{\sz_{\tilde i l}}{\wzd_l} \wf_l, 
$$
for all $0\leq i\leq m$.
Applying the morphism $F_*:\cf(\czcc)\ra \cfcc$ induced by the forgetful functor $F$ one obtains that
$$
\ch_i=\sum_{j=0}^{r} \frac{\sz_{\tilde i j}}{\wzd_{j}}F_*(\wf_j)=\sum_{j=0}^m\frac{\sz_{\tilde i \sg(j)}}{\widehat{d}_{\sg(j)}}F_j
.$$
This implies that $\muj(\ch_t)=\frac{\sz_{s \sg(j)}}{d_{\sg(j)}}=\wmu_{\sgj}(\wch_{\tilde t}),$ thus $\wmu_{\sgj}\in \ca_j$.
\epf

\newpage
\blue{Explain why $\mathbb Q\subseteq \mathbb K$ is a Galois extension? If it is... Check Harrison the second article.}

\blue{
\bne
\item
Define $\wcfcc$.  
\item
Define algebra structure on $\wcfcc$ and $\wpp_k(\jo,\jtw)$.\item 
Recall the isomorphism $\al:\wcfcc\ra \cecc$, right reference scbf8. It is defined in the Galois section.
\ene
}

\blue{ It follows that:
\bl
If $\cc$ is a weakly integral fusion category then
$$
m_\cc(\sgf(\ch), \sgf(\mu))=m_{\cc}(\ch, \mu)
$$
for any two class functions $\ch, \mu\in \cfcc$.
\el
}

\blue{
\bl
$$
\mu_j(\sg_F(\ch))=\sg(\mu_j(\ch))
$$
for any class function $\ch$.
\el
}

\subsection{Properties of the matrix $S$ in the modular case} Recall that the $S$-matrix of a braided pivotal fusion category is defined as $S_{XY}:=Tr(C_{X,Y}C_{Y,X})$. Then if $S_{ij}:=S_{M_i, M_j}$ it folllows
$$S_{ij}=S_{ji},\;\; S_{i^*j}=\ovl{S_{ij}}=S_{ij^*}$$
for all $i,j$.

As usually, by $\mathbb Q(S)$ is denoted  the field obtained by adjoining the $S$-matrix entries $S_{ij}$ to $\mathbb Q$. By \cite[Theorem 8.14.7]{EGNO15} one has that  $\mathbb Q(S)$ is contained in a cyclotomic extension.

By \cite[??]{ccc} in the modular case one has that
$\dim(\cc^k)=d_k^2$ and $\al_{ij}=\frac{S_{ij}}{d_j}$. The first orthogonality from Equation \eqref{611} can be written in this case as
\beq\label{611p}
\sum_{ k=0}^m S_{ik} S_{m^*k}=\delta_{i,m}\dimcc.
\eeq
\subsection{Comparison of the fields}
Since $S_{0j}=d_j$ one has that
$$
\mathbb Q\subseteq \mathbb K=\mathbb Q(\frac{S_{ij}}{d_j})\subseteq \mathbb Q(S)=\mathbb Q(S_{ij}).
$$
By the fundamental Theorem of Galois theory, since the extension $\mathbb Q\subseteq G_\cc=\mathbb Q(\frac{S_{ij}}{d_j})$ is Galois it follows that there is an epimorphism 
$$
G\xra{\pi} G_\cc,\;\; g\sent g|_{\mathbb K}.
$$
From here it follows that the two actions on $\mu_j$'s are compatible with $\pi$.

Denote also
\beq\label{f2}
s_{ij}:=\frac{S_{ij}}{\sqrt{\dimcc}}
\eeq
Then the first orthogonality relation can be written as
\beq\label{611ps}
\sum_{ k=0}^m s_{ik} s_{m^*k}=\delta_{i,m}.
\eeq

\subsection{On Theorem 8.14.7} Let $\cc$ be a modular tensor category. We use the following setup  $F_i:=\phir^{-1}(E_i)$ where $E_i$ is the primitive central idempotent of $X_i$ and $\phir$ is the Drinfeld map. By \cite[Theorem 8.13.11]{EGNO15} one has that
$$
\mu_i=\mu_{X_i},\;\; [Y]\sent \frac{S_{X_iY}}{d_i}.
$$
are all the morphisms $\mu:K_0(\cc)\ra \kk$.

With these notations one has
$$
\sg(\mu_{X_j})=\mu_{X_{{\tau(j)}}}\implies \sg(X_j)=X_{\tau(j)}$$
\blue{
The permutations on the simple objects from EGNO15 is my $\tau$.}

With the notations from the previous sections it follows that in this case:
$$
\al_{ij}=\mu_j([X_i])=\frac{S_{ij}}{d_j}=\frac{s_{ij}}{s_{0j}}.
$$
Then
$
\sg(\al_{ij})=\sg(\mu_j([X_i]))=\sg(\frac{S_{ij}}{d_j})=\frac{\sg(S_{ij})}{\sg(d_j)}.
$
On the other hand, by the definition of $\tau:=\tau_\sg$ one has
$
\sg(\al_{ij})=\al_{i\tau(j)}=\frac{S_{i\tau(j)}}{d_{\tau(j)}}.
$
Therefore
\beq\label{sgs}
\frac{\sg(S_{ij})}{\sg(d_j)}=\frac{S_{i\tau(j)}}{d_{\tau(j)}}.
\eeq
\blue{In terms of $s$-matrix this can be written
\beq\label{sgsp}
\frac{\sg(s_{ij})}{s_{0j}}=\frac{s_{i\tau(j)}}{s_{0\tau(j)}}.
\eeq
}
Switching $i$ and $j$ one has
\beq\label{sgs}
\frac{\sg(S_{ji})}{\sg(d_i)}=\frac{S_{j\tau(i)}}{d_{\tau(i)}}.
\eeq
Since $S_{ij}=S_{ji}$ the above two equations give that
\beq\label{f3}
\frac{S_{i\tau(j)}\sg(d_j)}{d_{\tau(j)}}=\frac{S_{j\tau(i)}\sg(d_i)}{d_{\tau(i)}}.
\eeq
Applying $\sg$ to the orthogonality relation (or directly apply Equation \eqref{sgccj})
$$
\sg(\frac{\dimcc}{d_i^2})=\sg(\frac{\dim(\cc)}{\dim(\cc^i)})\numeq{\ref{sgccj}}\frac{\dimcc}{\dim(\cc^{\tau(i)})}=\frac{\dimcc}{d_{\tau(i)}^2}
$$
which can be written as
\beq\label{sgccjmtc}
\sg\bigg(\frac{\dimcc}{d_i^2}\bigg)=\frac{\dimcc}{d_{\tau(i)}^2}.
\eeq
\bl
With the above notations one has $g(i^*)=g(i)^*$ and
\beq\label{gcls}
g(\frac{\dimcc}{d_i^2})=\frac{\dim(\cc)}{d_{g(i)}^2}
\eeq
\el
\bpf 
One has
\begin{eqnarray*}
g(\frac{\dimcc}{d_i^2})&=& g(\sumjtom \frac{S_{ij}}{d_i} \frac{S_{i^*j}}{d_{i^*}})
=\sumjtom g(\frac{S_{ij}}{d_i})g(\frac{S_{i^*j}}{d_{i^*}})=
 \\ &\numeq{\ref{sgs}}& \sumjtom \frac{S_{g(i)j}}{d_{g(i)}}\frac{S_{g(i^*)j}}{d_{g(i^*)}}\numeq{\ref{611p}} \delta_{g(i)^*, g(i)^*}\frac{\dim(\cc)}{d_{g(i)}^2}.
\end{eqnarray*}
This shows that $g(i^*)=g(i)^*$ and $g(\frac{\dimcc}{d_i^2})=\frac{\dim(\cc)}{d_{g(i)}^2}$.
\epf
Note that the above equation can be written in terms of $s$ as 
$$
g(s_{0i})^2=s_{0g(i)}^2.
$$
therefore there is a wekk defined sign such that 
\beq\label{defsgn}
g(s_{0i})=\eps_g(i)s_{0g(i)}
\eeq 
for all $i$.
\bl
One has 
\beq\label{tr2}
g(\frac{S_{ij}^2}{\dimcc})=\frac{S_{ig(j)}^2}{\dimcc}
\eeq
and
\beq\label{tr3}
g(s_{ij})=\eps_g(j)s_{ig(j)}=\eps_g(i)s_{jg(i)}.
\eeq
\el
\bpf
Indeed, 
\begin{eqnarray*}
g\bigg(\frac{S_{ij}^2}{\dim(\cc)}\bigg)&=&g\bigg((\frac{S_{ij}}{d_j})^2(\frac{d_j^2}{\dim(\cc)})\bigg)
\\&\numeq{\ref{sgs}, \ref{sgccjmtc}}&
(\frac{S_{i\tau(j)}}{d_{\tau(j)}})^2\frac{d_{\tau(j)^2}}{\dimcc}
\\&=&
\frac{S_{i\tau(j)}^2}{\dimcc}.
\end{eqnarray*}
This can be written as
\beq\label{sqij}
g\big(s_{ij}\big)^2=s_{i\tau(j)}^2.
\eeq
Using Equation \eqref{defsgn} it follows that
$$
g(s_{ij})=g\big(\frac{s_{ij}}{s_{0j}}s_{0j}\big)\numeq{\ref{sgsp}, \ref{defsgn}}\eps_g(j)\frac{s_{ig(j)}}{s_{0g(j)}}s_{0g(j)}=\eps_g(j)s_{ig(j)}.
$$
\epf

Moreover, the kernel of $\pi$ is $Gal(\mathbb Q(S)/Q)$ which coincides to all $g \in G$ such that $g(\frac{s_{ij}}{d_j})=\frac{s_{ij}}{d_j}$. Since by Equation \eqref{sgs} one has
$g(\frac{s_{ij}}{d_j})=\frac{s_{ig(j)}}{d_{g(j)}}$ it follows that the kernel coincides to all those elements that induces the identity permutation $\tau_{\pi(g)}$ on the set of all $\mu_j$'s.
\subsection{The action induced on $\cecc$}
The action induced on $\widehat{\cfcc}$ is translated by the  isomorphism $\al$ of algebras from \cite[Theorem 3.8]{scbf8} on $\cecc$. It follows that the action induced by $g \in G$ is the same as the action of $\pi(g)$.
\subsection{Self duality in the modular setting} By \cite[Theorem 4.1]{scbf8} it follows that the Drinfeld map $f_Q:\cfcc \ra \cecc$ is an isomorphism of algebras since it sends $\ch_i\sent \frac{\mtr C_i}{|\cc^i|}$.

\subsection{On the structure constants}
By \cite[Proposition 6.13]{ccc} the structure constants in the modular settings have the following formulae:
$$
c^k_{ij}=\frac{d_id_j}{d_k}N^k_{ij}
$$
\blue{Write also the formulae for $\wpp_k(\jo,\jtw)$.} Then
$$
\frac{d_id_j}{d_k}N^k_{ij}=c^k_{ij}=\frac{\dim(\cc^i)\dim(\cc^j)}{\dim(\cc^k)}{\wdht p}_k(i,j)=\frac{d_i^2d_j^2}{d_k^2}{\wdht p}_k(i,j).
$$
This gives that
\beq\label{wpkmtc}
{\wdht p}_k(i,j)=\frac{d_k}{d_id_j}N^k_{ij}
\eeq
\bp
If $\cc$ is a modular integral fusion category then $\wsgf$ is an algebra map. 
\ep
\bpf
By Equation \eqref{wpkmtc} it follows that in the integral case ${\wdht p}_k(i,j)\in \mathbb Q$.
Then the result  follows from item 3) of Proposition \ref{wsgfalgm}.
\epf
\br
In gnn it is applied the permutation of simple objects obtained from the permutation of the characters $\muj$.

In gnn one has
$$
\sg(S_{X Y})=S_{X, \sg(Y)}.
$$
\er
\br
Eventually work with $g \in \galkq$ instead of $\sg$.
\er
\subsection{rest of MTC}
On the other hand the second relation from Equation \eqref{612} becomes 
\beq\label{612p}
\sum_{ i=0}^mS_{il}S_{i^*k}=\delta_{l,k}\dimcc.
\eeq
which based on the above properties of the matrix $S$ is the same as the first relation of orthogonality.

\subsection{Questions in the modular settings} 

Q0: 
Are the two fields the same so we can speak of the same thing.

Q1:
Why the condition for $\sge$ algebra map implies the condition for $\wsgf$ an algebra map? At least in the modular case!

Q2: In the modular case does the permutation $\tau=\tau_\sg=\tau_\cc$ coincide with the permutation $\eta$ from my theory?

\section{To do}
To do:
\bne
\item Check that $G$ is abelian. Does it follow from the fact that is the Galois group of a sub-extension of the cyclotomic field?
\item Why the extension $\mathbb Q\subseteq \mathbb K$ is Galois?
\item Does the equality of multiplicities for any two classes follows by linearity? it does not make sense $\sg(\mu_j(\ch))$ for any class function $\ch$.
\item Permutation of the conjugacy classes, rational characters and rational conjugacy classes. read Berkovich.
\ene
\blue{
They are called algebraically conjugate characters in Berkovich.
\\
Two Galois conjugate characters have the same kernels and the same centers, also the same fields of values.
\\
In the group case one has $\sg(\ch)(g)=\sg(\ch(1))=\ch(1)$ since they are integers
Also 
$$
|\sg(z)|^2=\sg(z)\overline(\sg(z))=\sg(z)\sg(\overline(z))=\sg(|z|^2)=\sg(|z|)^2.
$$
\\ define rational characters and rational conjugacy classes.
}
\newpage
Write down the inverse of $\sg_F$ in terms of $\sg^{-1}_F$.
Indeed 
$$
\sg_F(F_j)=F_{\tau^{-1}_{\sg}(j)}.
$$ 
Then
$$
\sg^{-1}_F(\sg_F(F_j))=\sg^{-1}_F(F_{{\tau^{-1}_{\sg}}(j)})=F_{\tau^{-1}_{\sg^{-1}}({\tau^{-1}_{\sg}}(j)}=F_j
$$
\newpage
\section{Questions and comments}
\red{Galois theory for classes as in Bantay, The Galois permutation $\tau$ moves $g$-classes into $g$-classes.}

Q0: In defining the algebra structure on $\wcfcc$ one needs a fusion category in which the categorical dimensions are all $d_i\neq0$.

Q0' Need pivotal in order to speak about dimensions of simple objects. Otherwise one has to work with chosen isomorphisms.

Q1. Find an example where $\sge$ is an algebra map and $\wsgf(\ch_i)$ is not a pure character but a scalar of it. It should not be an integral category.

Q1' If it is weakly integral it follows that $d_i>0$ and therefore $d_i=d_{\eta(i)}$ by equation \eqref{}.

Q1: Ask Shimizu if one can get rid of pivotatity in Galois action working directly with the Grothendieck ring! More precisely, the isom

$$\al:\widehat{\grcc}\ra \cecc$$ is always an isomoprhism?

Q2': Quantum integral, i.e $d_i\in \mathbb Z$ implies integral in the usual sense? In the proposition with wsgf algebra map one needs quantum integral.

Q2: If $\sgf$ multiplicative, and bijective then it is automatically unitary?

Q3: Decide if one needs non-degenerate fusion category $\dimcc\neq 0$ in order to have the orthogonality relations. Over the complex field this is automatically non-degenarate.

Q4: Are all these true when $\cc$ is not pivotal?

Q5: Can I define the dual multiplication on $\wcfcc$ with $\fp$ instead of $\dim$? Which one has the coefficient of $\mu_0$ in $\muj\mu_{j\dl}$ positive in order to talk about a weak probability group?

Q6: Work in the proposition of $\wcfcc$ algebra map below first with $\cc$ pseudo-unitary, all dimensions are positive and algebraic integers. Can I then give up on this?

Q7: check that $d_i$ is an algebraic number and therefore integer if rational.

Q8: gnn proved the zero theorem without having integral quantum dimensions! So they don't need $\wsgf$ an algebra map.

\subsection{on the main theorem proof}

\blue{Check why the characters are permuted by $\sg$ in the weakly-integral modular settings, as claimed by GNN!}

\blue{Remove unnecessary lemmas and results concerning algebra maps conditions of $\sge$ and $\wsgf$!}

\blue{The proof is about $P_i$ which is a product indexed by $j$ to be rational. It is about the permutation $\tau$.}

Let $\cc$ be a integral fusion category with rational coefficients structure. It follows by Proposition \ref{wsgfalgm} that $\sge$ is an algebra map for all $\sg \in \galkq$. Since $\cc$ is also pseudo-unitary in this case it also follows by Corollary \ref{psu} that $d_i=d_{\eta(i)}$ for all $i$. Therefore, in this case Equation \eqref{sgebyev} becomes $\sgf(\ch_i)=\ch_{\eta(i)}$ for all $i$.

\br On the Proposition wsgf algebra map.

In other words, in vector matrix terms this can be written as
$$
\sg^{-1}(P)A=PAD
$$
where $(P)_{1k},\; (A)_{ki}=\al_{ik}\; D_{ii'}=\delta_{i, i'}\frac{d_i}{\sg^{-1}(d_i)}$. Since both $D$ and $A$ are invertible it follows we have equality.
\er
\newpage
\section{From the last section}
Let $\cc$ be a weakly integral fusion category.  By \cite[Proposition 8.27]{eno-annals} the dimensions of  simple objects in $\cc_{ad}$ are integers. Since $\cc^j$ are sum of simple object of the adjoint subcategory it follows that $\dimccj$ are integers.

It was shown in \cite[Proposition 8.24]{eno-annals} that there is a canonical spherical structure on $\cc$ such that categorical dimensions of objects in $\cc$ coincide with the Frobenius-Perron dimensions. We fix this spherical structure for the reminder of this section. Recall by \cite[Proposition 8.27]{eno-annals} that in such a category, the Frobenius–Perron dimension of any simple object is the square root of an integer.

As before, by $F_0, F_1, \dots ,F_m$ are denoted  the primitive central idempotents of $\cfcc$ and by $\mu_0, \mu_1, \dots \mu_m$  their corresponding characters on $\cfcc$. Therefore $\mui:\cfcc\ra \comp$ are algebras maps and $\mui(F_j)=\delta_{i,j}.$

We also denote by $\cc^0, \cc^1,\dots ,\cc^m$ the conjugacy classes of $\cc$ corresponding in this order to the primitive idempotents $F_0, F_1, \dots F_m$. Moreover, we let $\mtr C_0, \mtr C_1,\dots ,\mtr C_m$ denote their corresponding class sums.

Let also $M_0, M_1, \dots M_m$ be a complete set of representatives for the isomorphism classes of simple objects of $\cc$. As above, without loss of generality we may assume that $M_0=\unu$ is the unit object of $\cc$. Moreover, we denote by $\ch_i:=\mtr{ch}(M_i)$ the characters of the simple objects $M_i$.

Writing $\ch_i=\sumjtom\alij F_j$ we let $\mathbb K$ be the field extension of $\mathbf Q$ by all the scalars $\alij$. By \cite[Corollary  8.53]{eno-annals}  there is a cyclotomic field $\mathbb Q\subseteq Q(\xi)$ such that $\mathbb K\subseteq \mathbb Q(\xi)$.

Suppose that $\chio\chitw=\sum_{k=0}^mN^{k}_{\io \itw}\ch_{k}$. Expanding on the formula $\ch_{i}=\sumjtom\al_{ij}F_{j}$ it follows that
\beq\label{atfj}
\al_{\io j}\al_{\itw j}=\sum_{k}N^{k}_{\io\itw}\al_{k j}
\eeq

For any character $\muj:\cfcc\ra \comp$ we also define $\sg.\muj:\cfcc \ra \comp$  by
\beq\label{sgev}
[\sg.\mu_j](\ch_i):= \muj(\sg^{-1}(\ch_i))=\sg^{-1}(\alij)
\eeq

It is easy to see that 
 $\sg.\muj:\cfcc \ra \comp$ is an algebra map on ${\cfcc}$. Indeed one has that 
\begin{eqnarray*}
[\sg.\muj](\chio\chitw) &=&\sum_{k}N^{k}_{\io \itw}[\sg.\muj](\ch_{k})\numeq{\ref{sgev}}\sum_{k}N^{k}_{\io \itw}\sg^{-1}(\al_{kj})
\\ & =& \sg^{-1}(\sum_{k}N^{k}_{\io \itw}\al_{kj}) \numeq{\ref{atfj}} \sg^{-1}(\al_{\io j}\al_{\itw j})
\\ & = & \sg^{-1}(\al_{\io j})\sg^{-1}(\al_{\itw j})\numeq{\ref{sgev}} [\sg.\muj](\chio)[\sg.\muj](\chitw).
\end{eqnarray*}
\subsection{} We use the notation $\mtc J:=\{0,1,\dots ,m\}$.
\blue{Add the argument that $\tau$ is injective and therefore a permutation.}

Since $\sg.(\sg^{-1}.\muj)=\muj$ it follows that there is permutation $\tau=\tau_\sg$ on the set of indices $\mtc J$ such that $\sg.\muj:=\mu_{\tau (j)}$. It follows then from Equation \eqref{sgev} that for any $i,j$ one has:

\newpage

\br
Note that $\sg.\muj\neq \sg\circ \muj$ since 
$
(\sg.\muj)(\al\ch_i)=\al\sg(\muj(\ch_i))$ while 
$$
(\sg\circ \muj)(\al\ch_i)=\sg(\muj(\al\ch_i))=\sg(\al\muj(\ch_i))=\sg(\al)\sg(\muj(\ch_i)).
$$
\er

\br
The fact that $\sge$ is injective (and therefore bijective) can also be seen as follows: 
$$
\sge(z)=\sge(z')\implies <\ch_i, \sge(z)>=<\sgf(\ch_i), z>=<\ch_i, \sge(z')>=<\sgf(\ch_i), z'>$$
Therefore
$$<\sgf(\ch),z>=<\sgf(\ch),z'>$$ for all $\ch \in \cfcc$. Since the evaluation form is non-degenerate it follows that $z=z'$.
\er

\br
If $\sge$ is an algebra map then Equation \eqref{sgeqcls}
gives that $\eta(i^*)=\eta(i)^*$ and 
\beq\label{etaicond}
d_id_i^*=d_{\eta(i)}d_{\eta(i)^*},
\eeq
for all $i$.
\er

\red{
By Lemma \ref{eval2} the isomorphism $\al$ satisfies the following property
\beq\label{alev}
<\ch, \al(\mu)>=<\mu, \ch>.
\eeq
}

\red{\bt\label{main2} 
Suppose that $\cc$ is a weakly integral fusion category with rational structure coefficients. Then $\cc$ is integral.
\et
\bpf
By Lemma \label{ccwint} it follows that $\sge$ is an algebra map for all $\sg \in \galkq$.
\epf
}
\section{MTC-rest}

Moreover, in this case one has

$$
\al_{ij}=\mu_j([X_i])=\frac{S_{ij}}{d_j}.
$$
\subsection{On Theorem 8.14.7}
Morphisms
$$\mu_X:K_0(\cc)\ra \kk,\;\;Y\sent \frac{s_{XY}}{d_X}.$$
With my notations one has
$$\mu_{X_j}=\mu_{X_{{\tau(j)}}}\implies \sg(X_j)=X_{\tau(j)}$$

$$
\al_{ij}=\mu_j([X_i])=\frac{s_{ij}}{d_j}
$$
It follows that
$$
\sg(\al_{ij})=\sg(\mu_j([X_i]))=\sg(\frac{s_{ij}}{d_j})=\frac{\sg(s_{ij})}{\sg(d_j)}
$$
On the other hand 
$$
\sg(\al_{ij})=\al_{i\tau(j)}=\frac{s_{i\tau(j)}}{d_{\tau(j)}}
$$
It follows that
\beq\label{sgs}
\frac{\sg(s_{ij})}{\sg(d_j)}=\frac{s_{i\tau(j)}}{d_{\tau(j)}}
\eeq

Applying $\sg$ to the orthogonality relation (or directly Equation \eqref{sgccj})
$$
\sg(\frac{\dimcc}{\dim(X_i)^2})=\sg(\frac{\dim(\cc)}{\dim(\cc_i)})=\frac{\dimcc}{\dim(\cc_{\tau^{-1}(i)})}=\frac{\dimcc}{\dim(X_{\tau^{-1}(i)})^2}
$$
which can be written as
\beq\label{sgccjmtc2}
\sg(\frac{\dimcc}{d_i^2})=\frac{\dimcc}{d_{\tau^{-1}(i)^2}}.
\eeq
\blue{Apply the new trick:
$$
\sg(\frac{s_{ij}^2}{\dim(\cc)})=\sg((\frac{s_{ij}}{d_j})^2\frac{d_j^2}{\dim(\cc)})\numeq{\ref{sgs}, \ref{sgccjmtc}}(\frac{s_{i\tau(j)}}{d_{\tau(j)}})^2\frac{d^2_{{\tau^{-1}(j)}}}{\dimcc}
$$
}
\ed